\documentclass{amsart}

\usepackage{amsmath,amsthm,amsfonts,amssymb,bm,graphicx, mathabx,mathrsfs, dsfont}

\usepackage{comment}
\usepackage
{hyperref}

\hypersetup{colorlinks=true,citecolor=blue,linkcolor=blue,urlcolor=blue,
pdfstartview=FitH }

\theoremstyle{plain}
\newtheorem{lemma}{Lemma}
\newtheorem{theorem}[lemma]{Theorem}

\newtheorem{cor}[lemma]{Corollary}
\newtheorem{prop}[lemma]{Proposition}

\numberwithin{lemma}{section}
\numberwithin{equation}{section}

\theoremstyle{definition}

\theoremstyle{remark}

\begin{document}
\author{Edgar Assing}
\author{Valentin Blomer}
\author{Junxian Li}
 
\address{Mathematisches Institut, Endenicher Allee 60, 53115 Bonn}
\email{assing@math.uni-bonn.de}
\email{blomer@math.uni-bonn.de}

\address{Max-Planck-Institute f\"ur Mathematik,Vivatsgasse 7, 53111 Bonn }

\email{jli135@mpim-bonn.mpg.de}
  
\title{Uniform Titchmarsh divisor problems}

\thanks{The second author was supported in part by the DFG-SNF lead agency program grant BL 915/2-2}

\begin{abstract} Asymptotic formulae for Titchmarsh-type divisor sums are obtained with strong error terms that are uniform in the shift parameter. This applies to more general arithmetic functions such as sums of two squares, improving the error term in the representation of the number as a sum of a prime and two squares, and to Fourier coefficients of cusp forms, generalizing a result of Pitt. 
 
\end{abstract}


\setcounter{tocdepth}{2}  \maketitle 

\maketitle

\section{The Titchmarsh divisor problem}

The original divisor problem of Titchmarsh \cite{Ti} asks for the asymptotic evaluation of
\begin{equation}\label{1}
\sum_{a < p \leq x } \tau(p-a), \quad x \rightarrow \infty
\end{equation}
for a (positive) integer $a$, where $\tau$ denotes the standard divisor function and $p$ generally denotes a prime. For the past 90 years this problem has been a touchstone for the available techniques of analytic number theory. The obvious first step is to open the divisor function, which by Dirichlet's hyperbola method leads to the problem of counting primes in arithmetic progressions up to $\sqrt{x}$. Assuming the generalized Riemann Hypothesis (GRH) for all Dirichlet $L$-functions, Titchmarsh evaluated \eqref{1} as
\begin{equation}\label{2}
x\frac{\phi(a)}{a} \prod_{p \nmid a} \Big(1 + \frac{1}{p(p-1)}\Big)  + O_a\Big(\frac{x \log\log x}{\log x}\Big).
\end{equation}
Linnik \cite{Li} was the first to make this result unconditional using his dispersion method, but a much simpler route became possible  once the Bombieri-Vinogradov theorem was available, serving as a substitute for GRH for ``almost all'' Dirichlet $L$-functions. Short proofs for the equality of \eqref{1} and \eqref{2} along these lines can be found in \cite{Ro, Ha}. 

The next layer of  methodological depth came with the introduction of spectral methods of automorphic forms into analytic number theory  which allows us in certain situations to exhibit equidistribution of primes in arithmetic progressions to moduli slightly beyond $x^{1/2}$ via cancellation in sums of Kloosterman sums. Using combinatorial identities that go back to Linnik, Vaughan and Heath-Brown \cite{Li, Va, HB}, one can decompose the characteristic function on the primes into divisor-like functions. This leads to  weighted divisor sums in arithmetic progressions and after further transformations (Poisson summation) to sums of Kloosterman sums. 
Building on the ground-breaking work of Deshouillers-Iwaniec \cite{DI}, Fouvry \cite{Fo} and Bombieri-Friedlander-Iwaniec \cite{BFI} independently improved  the error term in \eqref{2} (with a slightly more precise main term) to
\begin{equation}\label{3}
O_{a, A}\big(x (\log x)^{-A}\big)
\end{equation}
  for an arbitrary number $A > 0$. In a sense this is the best we can hope for. Drappeau  \cite{Dr} showed  that any better error term would immediately improve our knowledge on Siegel zeros. Assuming GRH, he obtained a power saving error term
  \begin{equation}\label{4}
O_{a}\big(x^{1-\delta})
\end{equation}
for some $\delta > 0$. 

The Titchmarsh divisor problem is interesting and makes sense for many other arithmetic functions replacing the classical divisor function $\tau$. From the point of view of automorphic forms for the group ${\rm GL}(2)$, the most natural analogue of the divisor function are Hecke eigenvalues of cusp forms. It came as a surprise when Pitt \cite{P2} showed \emph{unconditionally}
\begin{equation}\label{5}
\sum_{p \leq x } \lambda(p-a)  \ll x^{1-\delta}
\end{equation}
for   $a = 1$, some $\delta > 0$ and Hecke eigenvalues $\lambda(n)$ of a holomorphic cuspidal Hecke eigenform for the group ${\rm SL}(2, \Bbb{Z})$. Although obviously different in some respects, the proof shares some features with the classical Titchmarsh divisor problem, in particular a rather similar appearance of sums of Kloosterman sums. \\

The strong error terms in \eqref{3}, \eqref{4}, \eqref{5} come with a price. While the equality of \eqref{1} and \eqref{2} can be shown in a large range of uniformity with respect to $a$  (although to our knowledge only the range $a \ll x^{1/2}$ seems to have appeared in print \cite{Ro}), all known techniques based on the Kuznetsov formula have a very restricted range for $a$. Drappeau obtains $a \ll x^{\delta}$ for some small, unspecified $\delta > 0$. To our knowledge the record is due to  Friedlander and Granville  \cite{FG}  who obtain the error term \eqref{3} uniformly in $a\ll x^{1/4 - \varepsilon}$ provided that $a$ has at most $O(\log\log x)$ prime factors (see also \cite{Fi}).  The uniformity restrictions are ultimately based on the same well-known restrictions in the distribution of primes in arithmetic progressions to moduli beyond $x^{1/2}$, cf.\ e.g.\ the discussion in \cite[p.\ 718]{Dr} and \cite{FG}. 

Also Pitt's method for the cuspidal version \eqref{5} of the Titchmarsh divisor problem is seriously restricted to small shifts $a$, as the discussion after \cite[(1.16)]{P2} shows. The reason is similar: in both the classical and the cuspidal case an important role is played by certain bilinear forms in Fourier coefficients of automorphic forms whose length and shape depends on the shift parameter $a$, and the analysis breaks down as soon as $a$ grows too quickly with $x$. The insufficient uniformity in $a$ implies in particular that the asymptotic evaluation of the equally natural dual sum
$$\sum_{p < N} \tau(N-p)$$
with an error term  of the quality \eqref{3} has been  out of reach, and the same applies to its cuspidal analogue. 

The present paper closes this gap. We start with a completely uniform Titchmarsh divisor formula. As usual, let $\Lambda$ denote the von Mangoldt function and let $\gamma$ be the Euler-Mascheroni constant. Let
\begin{equation}\label{c1c2}
\begin{split}
c_s(a) = \prod_{p \mid a} \Big(1 - \frac{1}{p^{s+1}}\Big) \prod_{p \nmid a} \Big( 1 + \frac{1}{(p-1)p^{s+1}}\Big).
\end{split}
\end{equation}

\begin{theorem}\label{thm-a} Let $\sigma \in \{\pm 1\}$,  $f\in \Bbb{Z}\setminus \{0\}$, $\max(-f, 0) \leq Y < X$ such that $X+f \asymp X$. If $\sigma =-1$, suppose that $f > X$. Then 
\begin{displaymath}
\begin{split}
&\sum_{Y < n \leq X} \Lambda(n) \tau(\sigma n + f)\\
& = c_0(f) \Big((X+\sigma f)\log(\sigma X+f) - (Y+\sigma f)\log(\sigma Y+f)  + (2\gamma - 1)(X-Y)\Big) +  2 c'_0(f)   (X-Y)
\\
& + O_A(X (\log X)^{-A})
 \end{split}
 \end{displaymath}
for any $A > 0$. If we assume GRH, then the error term can be improved to $O(X^{1-\delta})$ for  some  $\delta > 0$. 
\end{theorem}

This includes the original Titchmarsh divisor problem with long positive or negative shifts ($\sigma = 1$) as well as the dual problem ($\sigma = -1$). The same analysis can be applied successfully for other arithmetic functions in place of the divisor function, in particular convolutions $\chi_1 \ast \chi_2$ of two Dirichlet characters, see Theorem \ref{hooley} below.   As an application, we consider the classical problem of writing an integer $n$ as a sum of a prime and two squares. Hooley \cite{Ho} famously solved this under GRH, and Linnik \cite{Li1} used his dispersion method to obtain the first unconditional proof. Again a simpler proof can be obtained by the Bombieri-Vinogradov theorem, cf.\ e.g.\ \cite{EH}. In all cases, the error term saves only a small fractional power of $\log n$, which has not been improved since the 1960s. The following result goes much further. Let $\chi_{-4}$ denote the non-principal Dirichlet character modulo 4. 

\begin{theorem}\label{thm-b} 
The number of representations of an integer $n$ in the form $$n=p+x^2+y^2$$ in primes $p$ and non-zero integers $x, y$ is 
$$
\pi \, {\rm Li}(n) \prod_{p}\Big(1+\frac{\chi_{-4}(p)}{p(p-1)}\Big)\prod_{p\mid n}\frac{(p-1)^2}{p^2-p+\chi_{-4}(p)}+O_A\big(n(\log n)^{-A}\big)
$$
for any $A > 0$, where ${\rm Li}(n) = \int_2^n dt/\log t$ is the logarithmic integral. 
\end{theorem}



We finally treat the automorphic analogue  of  Titchmarsh's divisor problem with complete uniformity in the shift parameter, including the dual problem. We fix a (holomorphic or Maa{\ss}) cuspidal newform $\phi$ for $\Gamma_0(N)$ with   Hecke eigenvalues $\lambda(n)$. 

\begin{theorem}\label{thm-c} There exists $\delta > 0$  with the following property. Let $\sigma \in \{\pm 1\}$,  $f\in \Bbb{Z}\setminus \{0\}$, $\max(-f, 0) \leq Y < X$, such that $X+f \asymp X$. If $\sigma =-1$, suppose that $f > X$.  Then 
$$\sum_{Y < n \leq X} \Lambda(n) \lambda(\sigma n + f) \ll_{\phi} X^{1-\delta}$$
where the dependence on the conductor of $\phi$ is polynomial. 
\end{theorem}

A smooth version is derived in Proposition \ref{main2} below, from which Theorem \ref{thm-c} follows easily. 

\section{The methods and additional results}

While the statement of the Titchmarsh divisor problem and its variations are essentially elementary, the proofs use the full force of available machinery. This includes in particular bounds for multi-dimensional exponential sums based on the Riemann hypothesis for varieties (not only curves) over finite fields,  the spectral theory of automorphic forms, a subtle analysis of test functions in the Kuznetsov formula and representation theoretic methods for local considerations of Fourier coefficients. In this section we state some auxiliary results of independent interest. 

We introduce two fairly common pieces of notation. We write
$a \mid b^{\infty}$ to mean that $a$ has only prime divisors of $b$ and we denote by $\theta \leq 7/64$ an admissible exponent for the Ramanujan-Petersson conjecture. 

\subsection{Primes in arithmetic progressions}

The classical Titchmarsh divisor problem is based on equidistribution of primes in arithmetic progressions. The following result of independent interest implies Theorems \ref{thm-a} and \ref{thm-b} by a relatively straightforward procedure. 

\begin{theorem}\label{PrimeAp}
There exists a some positive constant $\delta$ with the following property. Let   $x\geq 2$, $C > 0$, $c, d \in \Bbb{N}$, $d \mid c^{\infty}$, $c_0,d_0 \in \Bbb{Z}$, $(d_0, d) = (c_0, c) = 1$, $a_1, a_2 \in \Bbb{Z} \setminus \{0\}$ such that
 $$Q\leq x^{1/2+\delta}, \quad a_1\leq x^{1 + \delta}, \quad a_2\leq x^\delta, \quad c, d \leq (\log x)^C$$  
 we have
 	\begin{align*}
	\sum_{\substack{q\leq Q\\ (q, a_1a_2)=1\\q\equiv c_0 \, (\text{{\rm mod  }}c)}}\Big(\sum_{\substack{n\leq x\\ n\equiv a_1\overline{a_2} \,(\text{{\rm mod  }}q)\\ n\equiv d_0(\text{{\rm mod  }}d)}}\Lambda(n)-\frac{1}{\phi(qd)}\sum_{\substack{n\leq x \\(n, qd)=1}}\Lambda(n)\Big)\ll_{C, A} x(\log x)^{-A}.
	\end{align*} 
If we assume GRH, then there exists $\delta' > 0$ such that in the longer range $c, d\leq x^{\delta'}$ we can bound the left hand side by $x^{1-\delta'}$. 
\end{theorem}

This should be compared, for instance, with \cite[Theorem 9]{BFI}, \cite[Theorem 1]{FG}, \cite[Theorem 6.2]{Dr}. 
The key point here is that $Q$ can be chosen a little larger than $x^{1/2}$ and in particular that  $a_1$ can be chosen as big as $x$ (and even a little bigger!). It is also useful to observe the additional congruence conditions on $q$ and $n$ which allow us to replace the divisor function in the Titchmarsh problem by sums of two squares and more general convolution functions, as the following result shows. 

\begin{theorem}\label{hooley}
	Let $X, Y, f, \sigma$ be as in Theorem \ref{thm-a}. Let $A, C > 0$ and let $\chi_1 \, (\text{{\rm mod }} c_1)$   and  $\chi_2\, (\text{{\rm mod }} c_2)$  be two primitive Dirichlet characters, where $c_1 = 1$ or $c_2 = 1$ is allowed. If  $\chi_1\not=\chi_2$ then we have for $c_1, c_2 \leq (\log X)^C$  that
		\begin{align*}
&	\sum_{Y < n \leq X}\Lambda(n)(\chi_1* \chi_2)(\sigma n+f)\\
	&=\Big(\frac{L(1, \chi_1\overline{\chi_2})c(\chi_1\overline{\chi_2}, f)\mu(c_2)}{\phi(c_2)}+\frac{L(1, \overline{\chi_1}\chi_2)c(\overline{\chi_1}\chi_2, f)\mu(c_1)}{\phi(c_1)} \Big)(X-Y) +O_{A, C}(X(\log X)^{-A}) 
	\end{align*}
	with
	\begin{align}\label{c-chi}
	c(\chi, f)=\prod_{p\nmid f}\Big(1+\frac{\chi(p)}{p(p-1)}\Big)\prod_{p\mid f}\Big(1-\frac{\chi(p)}{p}\Big).
	\end{align} 
 If we assume GRH, then the error term can be replaced by $O(X^{1-\delta'})$ for some $\delta'>0$ in the larger range $c_i\leq X^{\delta'}$.
\end{theorem}

As mentioned in the introduction, Theorem \ref{PrimeAp} can be obtained by combinatorial identities and suitable bounds for sums of Kloosterman sums, which is the topic of the following subsection.

\subsection{Sums of Kloosterman sums}

The key input for the proof of Theorem \ref{PrimeAp} is the following result.

\begin{theorem}\label{thm-dr} Let $C, D, N, R, S \geq 1$, $a, q, c_0, d_0 \in \Bbb{N}$ with $(c_0d_0, q) = 1$. Let $b_{n, r, s}$ be a sequence supported inside $(0, N] \times (R, 2R] \times (S, 2S] \cap \Bbb{N}^3$.   Let $g : \Bbb{R}^5_{>0} \rightarrow \Bbb{C}$ be a smooth function with compact support in $(C, 2C] \times (D, 2D] \times \Bbb{R}_{>0}^3$ such that$$\frac{\partial^{\nu_1 + \ldots + \nu_5}}{\partial c^{\nu_1} \partial d^{\nu_2} \partial n^{\nu_3} \partial r^{\nu_4} \partial s^{\nu_5} } g(c, d, n, r, s) \ll _{{\bm \nu}} ( c^{-\nu_1} d^{-\nu_2} n^{-\nu_3} r^{-\nu_4} s^{-\nu_5})^{1-\varepsilon_0}$$
for all ${\bm \nu} \in \Bbb{N}_0^5$ and some small fixed $\varepsilon_0 > 0$. Then
\begin{displaymath}
\begin{split}
&\sum_{\substack{c, d, n, r, s\\ c\equiv c_0, d \equiv d_0 \, (\text{{\rm mod }} q)\\ (qrd, sc)=1}} b_{n, r, s} g(c, d, n, r, s) e\Big( an \frac{\overline{rd}}{sc}\Big) \ll_{\varepsilon, \varepsilon_0} (aqCDNRS)^{\varepsilon + O(\varepsilon_0)} q^{3/2} K(a, C, D, N, R, S) 
\end{split}
\end{displaymath}
where
\begin{displaymath}
\begin{split}
K(a, C, D, N, R, S)^2 = &a^{2\theta} q^{2}  \sum_{\substack{n'' \mid a^{\infty}\\ n'' \leq 2N}} (n'')^{2\theta}  \Big(CS\Big(RS+\frac{N}{n''}\Big)( C+RD)  + aSNR\Big) \| \tilde{\textbf{b}}(n'') \|_2^2 \\
&+ q  \big( C^2 D S \sqrt{R(RS + N)}   + D^2 N R\big) \| \textbf{b} \|_2^2 . 
\end{split}
\end{displaymath}
with  
$$\tilde{b}_{n, r, s}(n'') =  b_{nn'', r, s}.$$
\end{theorem}

Bounds of this kind have a long history starting with the ground-breaking work of Deshouillers-Iwaniec \cite[Theorem 12]{DI}. For $a = 1$, Theorem \ref{thm-dr} recovers \cite[Theorem 2.1]{Dr} except from unimportant powers of $q$ and a similarly unimportant inaccuracy in the last term $D^2NR\|\textbf{b} \|^2$ where in \cite[Theorem 2.1]{Dr} no factor $1/S$ should occur (cf.\ \cite{BFI3}). Note that the sequence  $ \tilde{\textbf{b}}(n'')$ has  typically fewer terms than $\textbf{b}$. In many practical cases, the factor $\| \tilde{\textbf{b}}(n'') \|_2^2$ should decrease with $n''$ and therefore amortize the $n''$-sum. The novelty of Theorem \ref{thm-dr} is the additional parameter $a$ which is not restricted in any way. This offers a lot of extra flexibility in applications. 

An experienced reader can guess the strategy of the proof from the shape of the bound. The Kuznetsov formula translates the quintilinear form in  Kloosterman fractions into a spectral expression involving sums 
$$\sum_{N < n \leq 2N} a_n \sqrt{an} \rho_{f, \mathfrak{a}} (an)$$
for a general   sequence $a_n$, a cusp form $f$ with Fourier coefficients $\rho_f$ and a cusp $\mathfrak{a}$ which will ultimately be $\mathfrak{a} = 1/s$.  In order to apply the spectral large sieve \cite{DI} in an efficient way, one would like to ``factor out'' $a$ from $\rho_{f, \mathfrak{a}} (an)$ which eventually leads to the extra factor $a^{2\theta}$  in the bound for $K(a, C, D, N, R, S)^2$. The first observation is that this sacrifice of a factor $a^{2\theta}$ with any value of $\theta < 1/8$ still suffices for the application to the Titchmarsh divisor problem. We highlight that this is one of the rare applications where a particularly strong approximation towards the Ramanujan-Petersson conjecture is needed. 
The second observation is that the desired factorization is, however, not easily possible, among other things because Fourier coefficients at arbitrary cusps are not necessarily multiplicative, certainly not completely multiplicative, and even if they are for newforms, one also has to take care of oldforms in the spectrum. The key input here is Lemma \ref{lem1-dr}  below.

For the proof of Theorem \ref{thm-c} we need a different configuration of Kloosterman sums. The following is  a  generalization of \cite[Theorem 1.6]{P2}.
\begin{theorem}\label{lem5} Let  $s, r, n \in \Bbb{N}$ with $(r, s) = 1$. Let $M, C, Z \geq 1$. Let $F$ be a function with support in $[M, 2M] \times [C, 2C]$ satisfying $F^{(\nu_1, \nu_2)} \ll_{\nu_1, \nu_2} Z^{\nu_1+\nu_2} M^{-\nu_1} C^{-\nu_2}$ for all $\nu_1, \nu_2 \in \Bbb{N}_0$. Write
$X = (Mn/(s^2r C^2))^{1/2}$ and suppose that $X \ll Z$. Let $\alpha_m$ be any sequence of complex numbers. Then
\begin{displaymath}
\begin{split}
&   \sum_{(c, r) = 1} \frac{1}{c} \sum_m \alpha_m F(m, c) S(m \bar{r}, \pm n, sc) \\
&\ll  \| \alpha \| Z\Big(Z (sM)^{1/2}+ Z^2 s\sqrt{r} + \frac{(M(n, rs))^{1/4}(ZsC)^{1/2} }{r^{1/4}} + \frac{Z^{3/2} (n, rs)^{1/4}C^{1/2} r^{1/4} s}{M^{1/4}}\Big) (CMZnrs)^{\varepsilon}.
\end{split}
\end{displaymath}
\end{theorem}
 The key point   is that unlike \cite[Theorem 1.6]{P2} we do not require the assumption $n \leq rs$. The proof leads to similar sums of Fourier coefficients, but here we cannot afford to lose a factor of size $n^{\theta}$, even for $\theta$ as small as $7/64$.   Instead, we apply the Kuznetsov formula backwards with a different choice of test function. The arising Kloosterman sums can then be bounded successfully by Weil's bound. Interestingly, this strategy would not succeed for the application to the classical Titchmarsh divisor problem in Theorem \ref{thm-a}.

\subsection{Automorphic shifted convolution problems} An important role in the proof of Theorem \ref{thm-c} is played by correlations of various additively shifted arithmetic functions. The application of Vaughan/Heath-Brown identities reduces the characteristic function of primes to generalized divisor functions. Pitt's proof of \eqref{5} is based on a power saving in type III sums. 

 Let $\phi$ be a  (holomorphic or Maa{\ss}) cuspidal newform  for $\Gamma_0(N)$ with   Hecke eigenvalues $\lambda(n)$. If $\phi$ is holomorphic we write $k$ for the weight and if it is a Maa{\ss} form we write $t \in \Bbb{R} \cup[0, 7i/64]$ for the spectral parameter. For uniform notation we write $\mu = k$ or $1 + |t|$ depending on the case.

The following is a generalization of \cite[Theorem 2]{P1} which is a smooth version of a shifted convolution sum of Hecke eigenvalues $\lambda(n)$ with the ternary divisor function $\tau_3$. As before, this result is completely uniform in all considered parameters. The smoothing can easily be removed if desired.

\begin{theorem}\label{thm1} Let $r, f\in \Bbb{Z} \setminus \{0\}$, $Z, X_1, X_2, X_3 \geq 1/100$. 
Suppose that at least one of $r, f$ is positive. Let $G$ be a function with support on  $\bigtimes [X_j, 2X_j]$ and $G^{(\nu_1, \nu_2, \nu_3)} \ll_{{\bm \nu}} Z^{\nu_1+\nu_2+\nu_3}X_1^{-\nu_1} X_2^{-\nu_2} X_3^{-\nu_3}$ for all $\bm \nu \in \Bbb{N}_0^3$.  Write $X =  |r|X_1X_2X_3$ and suppose that $G(k, l, m)$ vanishes unless $rklm + f$ is in some dyadic interval $[\Xi, 2\Xi]$ with $X/Z \ll \Xi \ll X+|f|$. Then 
\begin{equation*}
\begin{split}
\sum_{k, l, m} G(k, l, m)\lambda(rklm + f) \ll  (\mu Z | f| NX)^{\varepsilon}  Z^3N^{5/2}\mu^{7/2}  (1 + |f|/X)^{25/24}  (X + |f|)^{23/24} .
\end{split}
\end{equation*}
\end{theorem}

We turn to the shifted convolution problem of two Hecke eigenvalues and consider a situation similar to \cite[Theorem 1.4]{P2}. In the special case $f=1$, this matches exactly \cite[Theorem 1.4]{P2} with explicit polynomial dependence on $\mu, Z, N$, but unlike \cite[Theorem 1.4]{P2} we allow shifts as long as $X$ (and even a little longer).  Therefore the proof will turn out to be quite different in several aspects. 

\begin{theorem}\label{lem8} Let $M \geq 1$ and let $m_1, m_2 \in \Bbb{Z}\setminus \{0\}$ be distinct, squarefree and of the same sign satisfying $M \leq |m_1|, |m_2| \leq 2M$. Let $f \in \Bbb{Z} \setminus \{0\}$ and suppose that $m_1, m_2, f$ are not all negative. Let $d = (m_1, m_2)$, $m_1 = dm_3$, $m_2 = dm_4$ and $D = (m_1 - m_2, d^{\infty})$. Let $X, Z \geq 1$ and $g$ be a function with support in $[X, 2X]$ and $g^{(\nu)} \ll_{\nu}  (X/Z)^{-\nu}$ for all $\nu \in \Bbb{N}_0$. Suppose that $g(n)$ vanishes unless $m_1n+f, m_2n+f$ are in some dyadic interval $[\Xi, 2\Xi]$ with $MX/Z \ll\Xi \ll MX + |f|$. 
Then
\begin{displaymath}
\begin{split}
\sum_n \lambda(m_1n + f) \lambda(m_2n+f) & g(n) \ll ( |f| X M Z N \mu)^{\varepsilon} \mu^{23/2}  \Big (\mu + Z \Big(1 + \frac{|f|}{MX}\Big)\Big)^8 N^{5/2} \\
&\Big( \frac{(MX + |f|)^{1/2}M^{1/2} D}{d^{1/2}} + \frac{(MX+|f|)^{3/4}D^{1/4} (f, Dm_3m_4)^{1/4} }{M^{1/4}}\Big). 
\end{split}
\end{displaymath}
\end{theorem}

The technology developed in this paper has  applications to other arithmetic problems. Drappeau \cite[Theorem 1.5]{Dr} obtained an asymptotic formula with a power saving error term for the shifted convolution problem
$$\sum_{n \leq x}\tau_k(n)\tau(n+h)$$
for $h=1$ and remarks that   shifts of length $h \ll x^{\delta}$ for some $\delta > 0$ are also possible. Topacogullari \cite[Theorem 1.1]{To} obtains an explicit result for $h \ll x^{15/19}$, and it is easy to see from \cite[(2.2) - (2.4)]{To} that his result can be upgraded to $h \ll x^{1-\varepsilon}$. The methods of this paper, in particular Proposition \ref{TypeII} (which is a consequence of Theorem \ref{thm-dr}), imply a result where $h$ can be as large as $x$ and even a little bigger. The details will appear elsewhere. \\

\textbf{Roadmap for the rest  of the paper:} In Section \ref{sec3} we prove Theorem \ref{thm-dr}. In Section \ref{sec4} we deduce Theorem \ref{PrimeAp} from Theorem \ref{thm-dr}. In order to avoid undue redundancy, we assume some familiarity with the paper of Drappeau \cite{Dr} in these two sections. It is in principle simple to obtain Theorems \ref{thm-a}, \ref{hooley} and \ref{thm-b} from Theorem \ref{PrimeAp}, but the deduction of the respective main terms is technically a bit challenging.  The relevant computations  are the content of Section \ref{sec5}. Next we prepare the scene for the proof of Theorem \ref{thm-c}, where we assume some familiarity with the papers of Pitt \cite{P1, P2}.  In Section \ref{sec6} we prove Theorem \ref{lem5}. This is followed by a number of preparatory lemmas in Section \ref{sec7} that we need in order to prove Theorems \ref{thm1} and \ref{lem8} as well as some consequences in Section \ref{sec8}. At this point we have all ingredients available to complete the proof of Theorem \ref{thm-c} in Section \ref{sec9}.

\section{Proof of Theorem \ref{thm-dr}}\label{sec3}

\subsection{A reduction step} We start with the following argument which allows us to assume that   $b_{n, r, s}$ vanishes unless $(a, rs) = 1$. Indeed, writing $b_1 = (a, r)$, $b_2 = (a, s)$, $a = b_1 b_2 b_3$, $r = r'b_1$, $s = s'b_2$, the sum in question equals
\begin{displaymath}
\begin{split}
&\sum_{ \substack{b_1b_2b_3 = a \\ (b_2, qb_1) = 1}} \sum_{\substack{c, d, n, r', s'\\ c\equiv c_0, d \equiv d_0 \, (\text{{\rm mod }} q)\\ (qr' d, s' c) = 1\\ (r', b_2b_3) = (s', b_1b_3) = 1\\ (c, b_1) = (d, b_2) = 1}} b_{n, r'b_1, s'b_2} g(c, d, n, r'b_1, s'b_2) e\Big( b_3n \frac{\overline{r'd}}{s'c}\Big)\\
=&\sum_{ \substack{b_1b_2b_3 = a \\ (b_2, qb_1) = 1}} \sum_{\substack{f_1 \mid b_1\\ f_2 \mid b_2\\ (f_1f_2, q) = 1}} \mu(f_1)\mu(f_2) \sum_{\substack{c, d, n, r', s'\\ c\equiv \bar{f}_1c_0, d \equiv \bar{f}_2d_0 \, (\text{{\rm mod }} q)\\ (qr' f_2d, s' f_1c) = 1\\ (r', b_2b_3) = (s', b_1b_3) = 1}} b_{n, r'b_1, s'b_2} g(cf_1, df_2, n, r'b_1, s'b_2) e\Big( b_3n \frac{\overline{r'f_2d}}{s'f_1c}\Big).
\end{split}
\end{displaymath}
We have 
 $$ 
 e\Big( b_3n \frac{\overline{r'f_2d}}{s'f_1c}\Big) =  e\Bigg( \frac{b_3}{(b_3, f_1f_2)} n \frac{\overline{r'\frac{f_2}{(f_2, b_3)}d}}{s'\frac{f_1}{(f_1, b_3)}c}\Bigg).$$
Now  the   result for general sequences $b_{n, r, s}$ follows from the special result for sequences with $(a, rs) = 1$ with $$\Big(\frac{C}{f_1}, \frac{D}{f_2}, \frac{Rf_2}{(f_2, b_3)b_1}, \frac{Sf_1}{(f_1, b_3)b_2}, \frac{b_3}{(b_3, f_1f_2)}\Big)$$ in place of $(C, D, R, S, a)$ and the sequence $$b'_{n, r, s} = \begin{cases}b_{n, r'b_1, s'b_2}, & r = r'\frac{f_2}{(f_2, b_3)}, s = s'\frac{f_1}{(f_1, b_3)},    (r', b_2b_3) = (s', b_1b_3) = 1,\\ 0, & \text{else,}\end{cases} $$ 
where we have
$$\Big(rs,\frac{b_3}{(b_3, f_1f_2)}\Big) = \Big( \frac{r's'f_1f_2}{(f_1f_2, b_3)}, \frac{b_3}{(b_3, f_1f_2)}\Big) = 1$$
on the support of $b'_{n, r,s}$.  

\subsection{Spectral inequalities} From now on we assume that $b_{n, r, s} = 0$ unless $(rs, a) = 1$. We write 
\begin{equation}\label{Q}
Q = q r s = Q'Q'', \quad (s, rq) = 1, \quad Q' = (Q, (rs)^{\infty}), \quad (Q'', rs) = 1.
\end{equation}
In particular, $(a, Q') = 1$. 
We follow the proof of Drappeau \cite[Theorem 2.1]{Dr}, which after Poisson summation in $d$ reduces the problem to bounding sums over Kloosterman sums (in arithmetic progressions). We use the same notation and normalization as Drappeau and define corresponding to \cite[(4.23) -- (4.25)]{Dr} the three spectral quantities
\begin{equation}\label{spectral}
\begin{split}
&  \mathcal{H}_{\mathfrak{a}}(T, Q, \chi, N, (a_n), a) = \sum_{\substack{\kappa < k \leq T\\ k \equiv \kappa \, (\text{mod } 2)}} \Gamma(k) \sum_{f \in \mathcal{B}_k(Q, \chi)} \Big| \sum_{N < n \leq 2N} a_n \sqrt{an} \rho_{f , \mathfrak{a}}(an)\Big|^2,\\
&  \mathcal{M}_{\mathfrak{a}}(T, Q, \chi, N, (a_n), a) =  \sum_{\substack{f \in \mathcal{B}(Q, \chi)\\ |t_f| \leq T}}\frac{(1 + |t_f|)^{\pm \kappa}}{\cosh(\pi t_f)} \Big| \sum_{N < n \leq 2N} a_n \sqrt{an} \rho_{f , \mathfrak{a}}(\pm an)\Big|^2,\\
&  \mathcal{E}_{\mathfrak{a}}(T, Q, \chi, N, (a_n), a) = \sum_{\mathfrak{c} \text{ singular}} \int_{-T}^T \frac{(1 + |t_f|)^{\pm \kappa}}{\cosh(\pi t_f)} \Big| \sum_{N < n \leq 2N} a_n \sqrt{an} \rho_{\mathfrak{c} , \mathfrak{a}}(\pm an, t)\Big|^2
\end{split} 
\end{equation}
for a singular cusp $\mathfrak{a}$, $T, N \geq 1$, $a, q\in \Bbb{N}$, $\chi$ a Dirichlet character modulo $q_0 \mid q$ and parity $\kappa$,   $a_n$ a sequence of complex numbers, and $\mathfrak{c}$ runs over singular cusps for $\Gamma_0(Q)$ with respect to $\chi$. We only need the two cases $\mathfrak{a} \in \{\infty, 1/s\}$, and we will have $a=1$ if $\mathfrak{a} = \infty$. The complicated part is to treat the case $\mathfrak{a} = 1/s$ for general $a$.   

The basic idea is now to assume, as we may, that $\mathcal{B}_k(q, \chi), \mathcal{B}(q, \chi)$ are Hecke eigenbases and to factor out $a$ from $\rho_{f, 1/s}(an)$, and similarly for Eisenstein series. This comes with considerable difficulties. First of all, even for a newform $f$ the Fourier coefficients at the cusp $1/s$ are not multiplicative in a strict sense. Secondly, the orthonormal bases also contain oldforms, and common factors of $a$ with the level of $f$ cause problems. 

We start with a detailed discussion of Maa{\ss} forms, the case of holomorphic forms and Eisenstein series is slightly easier and requires only minor notational modifications. We denote by $\mathcal{B}^{\ast}(Q_1, Q, \chi)$   an orthonormal basis of \emph{newforms} of level $Q_1$ with $q_0 \mid Q_1 \mid Q$ normalized with respect to $\Gamma_0(Q)$. For any Maa{\ss} form $f$ and $d \in \Bbb{N}$ we write $f_d(z) = f(dz)$. Note that this map is an isometry, i.e.\ $\| f_d \| = \|f \|$, if both inner products are taken with respect to the same group. For $f\in \mathcal{B}^{\ast}(Q_1, Q, \chi)$   we choose an orthonormal basis $\{f^{(g)} : g \mid Q/Q_1\}$ in the space $\langle \{f_d : d\mid Q/Q_1\} \rangle$ as in \cite[Section 5]{BM}. Then we have 
\begin{equation}\label{0new}
f^{(g)} = \sum_{d \mid g} \xi_g(d) f_d, \quad \xi_g(d) \ll (g/d)^{\theta - 1/2+\varepsilon} \ll 1. 
\end{equation}
Using this particular basis and the Cauchy-Schwarz inequality, we can re-write $ \mathcal{M}_{1/s}(T, Q, \chi, N, (a_n), a) $ as follows
\begin{equation}\label{3-line}
\begin{split}
&\sum_{q_0 \mid Q_1 \mid Q}  \sum_{g \mid \frac{Q}{Q_1}} \sum_{\substack{f \in \mathcal{B}^{\ast}(Q_1,Q,  \chi)\\ |t_f| \leq T}}\frac{(1 + |t_f|)^{\pm \kappa}}{\cosh(\pi t_f)} \Big| \sum_{N < n \leq 2N} a_n  \sum_{d \mid g} \xi_g(d) \sqrt{an}\rho_{f_d , 1/s}( \pm an)\Big|^2 \\
& \leq \sum_{q_0 \mid Q_1 \mid Q}  \sum_{g \mid \frac{Q}{Q_1}} \sum_{\substack{f \in \mathcal{B}^{\ast}(Q_1,Q,  \chi)\\ |t_f| \leq T}}\frac{(1 + |t_f|)^{\pm \kappa}}{\cosh(\pi t_f)}\sum_{d\mid g} |\xi_g(d)|^2 \sum_{d\mid g}\Big| \sum_{N < n \leq 2N} a_n   \sqrt{an}\rho_{f_d , 1/s}( \pm an)\Big|^2 \\
& \ll Q^{\varepsilon} \sum_{q_0 \mid Q_1 \mid Q}  \sum_{d\mid \frac{Q}{Q_1}} \sum_{\substack{f \in \mathcal{B}^{\ast}(Q_1, Q, \chi)\\ |t_f| \leq T}}\hspace{-0.2cm}\frac{(1 + |t_f|)^{\pm \kappa}}{\cosh(\pi t_f)}   \Big| \sum_{n'' \mid a^{\infty}}  \sum_{\substack{N/n'' < n' \leq 2N/n''\\ (n', a) = 1}} a_{n'n''}  \sqrt{an'n''}\rho_{f_{d} , 1/s}( \pm an'n'')\Big|^2.
\end{split}
\end{equation}
To proceed further, we need the following lemma. 
\begin{lemma} \label{lem1-dr} Let  $f\in \mathcal{B}^{\ast}(Q_1, Q, \chi)$, $d \mid Q/Q_1$,  $(Q/s, s) = 1$. Let $(n, m) = 1$ and define
$$K = \Big( \frac{(s, d)Q/s}{(dQ_1, Q/s)}, n^{\infty}\Big).$$
Then
$$\sqrt{nm}\rho_{f_d, 1/s}(\pm nm) = \delta_{K \mid n} \lambda\Big(\frac{n}{K}\Big) e\Big( \pm (nm - Km)\frac{\overline{s}}{Q/s}\Big)\sqrt{Km} \rho_{f_d, 1/s}(\pm Km)$$
where $\lambda$ is Hecke eigenvalue of a certain cusp form depending on $f$, $s$, $d$ and $Q$.  
\end{lemma}

\emph{Remarks:} An explicit formula for $\lambda$ is given in the proof. Somewhat less general results in this direction can be found in \cite{GHL}.\\

\textbf{Proof.} We start by setting up some notation. Write $M=Q/s$ and decompose $d=d_1d_2$ with $d_1=(d,s)$ and $d_2=d/d_1$. By assumption we have $(M,s)=(d_1,d_2)=1$. We also define $P_1=(d_2Q_1,s/d_1)$ and $P_2= d_2Q_1/P_1 $. We claim that $(P_1,P_2)=1$. Indeed, take $p\mid P_1$ then $p\mid s$ and the $p$-adic valuation of $s$ is given by $v_p(s) = v_p(Q)$. Since $d\mid  Q/Q_1$ and $(d_1,d_2)=1$, we have $v_p(d_2Q_1)\leq v_p(s/d_1)$. In particular, $v_p(P_1) = v_p(d_2Q_1)$ and $p\nmid P_2$ as claimed. Next, we claim that $d_2\mid P_2$. Instead of showing this directly we compute $(d_2,P_1)=1$, which amounts to the same. To do so take $p\mid d_2$. We observe that $v_p(P_1) = \min(v_p(d_2Q_1),v_p(s/d_1)) = 0$. The last equality follows since by construction $(d_2,s) = 1$. Finally we observe that $(d_1, P_2) = 1$ since $P_2 = (d_2Q_1, Q/s)$ and $(s, Q/s) = 1$. 
For later reference we summarize: 
\begin{equation}
	d_2\mid P_2\mid d_2Q_1 \quad \text{ and }\quad  \Big(P_2,\frac{d_2Q_1}{P_2}\Big) = (d_1, P_2) = 1.\label{eq:divisibility_prop}
\end{equation}
Finally, we set $L= (s,d)M/(dQ_1,M)$. It is easy to check that $L= d_1M/P_2$.

Since $\chi$ has conductor $q_0\mid Q_1$, we can view $\chi$ as a character modulo $d_2Q_1$ or $Q_1$ when convenient. Furthermore, there are Dirichlet characters $\chi_{1}$ with conductor $(P_1,q_0)$ and $\chi_{2}$ with conductor $(P_2,q_0)$ such that $\chi = \chi_{1}\chi_{2}$.

Pick $\alpha,\beta\in\mathbb{Z}$ such that 
\begin{equation}
	P_2\beta-\frac{s}{d_1}\alpha = 1. \nonumber
\end{equation}
This is possible because $P_2$ was constructed such that $(P_2,s/d_1)=1$. With this at hand, we decompose 
\begin{equation}
	\frac{1}{\sqrt{d_1M}}\left(\begin{matrix} d_1M & 0 \\ Q & 1 \end{matrix}\right) = \frac{1}{\sqrt{d_1M}}\left(\begin{matrix} 1 & \alpha \\ \frac{s}{d_1} & P_2\beta \end{matrix}\right)\left(\begin{matrix} d_1M & -\alpha \\ 0 & 1 \end{matrix}\right).\nonumber
\end{equation}
Let
\begin{equation}
	\sigma_{\mathfrak{a}} := \frac{1}{\sqrt{M}}\left( \begin{matrix}  M & 0 \\ Q &  1 \end{matrix}\right) \nonumber 
\end{equation}
be the scaling matrix for the cusp $\mathfrak{a}=1/s$ as in \cite[(4.1)]{Dr}. Using the decomposition above we find 
\begin{equation}
	f_d(\sigma_{\mathfrak{a}}z) = f_{d_2}\left(\left(\begin{matrix} 1 & \alpha \\ \frac{s}{d_1} & P_2\beta \end{matrix}\right).[d_1Mz-\alpha]  \right). \nonumber
\end{equation}

As in \cite[Lemma~A.3]{KMV}, we can write
\begin{equation}
	\left(\begin{matrix} 1 & \alpha \\ \frac{s}{d_1} & P_2\beta \end{matrix}\right) = \gamma W_{P_2} \left(\begin{matrix} P_2^{-1} & 0 \\ 0 & 1 \end{matrix}\right) \nonumber
\end{equation}
for $\gamma\in\Gamma_0(d_2Q_1)$ and a matrix
\begin{equation}
	W_{P_2} = \left(\begin{matrix} xP_2 & y \\ zd_2Q_1 & wP_2 \end{matrix}\right) \text{ with } P_2xw-P_1yz = 1 ,\,  y\equiv 1 \text{ mod } P_2, \, x\equiv 1\text{ mod } P_1, \nonumber
\end{equation}
representing the $P_2$-Atkin-Lehner operator. Unfolding this matrix decomposition  and considering it modulo $P_1$ and $P_2$ reveals $\chi(\gamma) = \chi_{2}(-s/d_1)$. We arrive at
\begin{align}
	f_d(\sigma_{\mathfrak{a}}z) &= \chi_{2}(-s)\overline{\chi_{2}(d_1)}f_{d_2}\left(W_{P_2}.\left[\frac{d_1M}{P_2}z-\frac{\alpha}{P_2}\right]  \right) \nonumber \\ 
	&= \chi_{2}(-s)\overline{\chi_{2}(d_1)\chi_{1}(d_2)}f\left(W_{P_2/d_2}.\left[Lz-\frac{\alpha}{P_2}\right]  \right). \nonumber
\end{align}
In the last step we used \cite[Proposition~1.5]{AL} (which of course remains true   for Maa\ss\  forms).  Note that the required conditions are satisfied by \eqref{eq:divisibility_prop}.

According to \cite[Proposition~A.1]{KMV} (see also \cite{AL}), there is an arithmetically normalized newform $g$ of level $Q_1$, nebentypus $\chi_1\overline{\chi_2}$ and the same spectral data such that
\begin{equation}
	f(W_{P_2/d_2}z) = \rho_{f, \infty}(1)\eta_f(P_2/d_2) g(z),\nonumber
\end{equation}
where $\eta_f(P_2/d_2)$ is the Atkin-Lehner pseudo-eigenvalue of $f$. Note that $g$ is uniquely determined by $f$ and $P_2/d_2$. We conclude that \begin{equation}
	f_d(\sigma_{\mathfrak{a}}z) =  C_{f,s,d}g\left(Lz-\frac{\alpha}{P_2} \right) \label{eq:A-L_identity}
\end{equation}
for $C_{f,s,d} = \chi_{2}(-s)\overline{\chi_{2}(d_1)\chi_{1}(d_2)}\rho_{f, \infty}(1)\eta_f(P_2/d_2)$. The Fourier coefficients $\rho_{g,\infty}(n)$ are multiplicative and directly relate to the Hecke eigenvalues of $g$, which in turn can be explicitly described in terms of the Hecke eigenvalues of $f$. In particular, $\rho_{g,\infty}$ can be estimated  using any available bound towards the Ramanujan conjecture for Maa\ss\  forms. 

Expanding both sides of \eqref{eq:A-L_identity} in the corresponding Fourier expansions and comparing coefficients yields the identity
\begin{equation}
	\rho_{f_d,1/s}(n) = \delta_{L\mid n}C_{f,s,d} e\left(-n\frac{\alpha}{LP_2}\right)\rho_{g,\infty}\left(\frac{n}{L}\right). \nonumber 
\end{equation}
The stated formula, with $\lambda(n) = \sqrt{n}\rho_{g,\infty}(n)$, follows immediately from this equation after observing that
\begin{equation}
	e\Big(-\frac{n}{L}\cdot\frac{\alpha}{P_2}\Big) = e\Big(n\frac{\overline{s/d_1}}{d_1M}\Big) = e\left(n\frac{\overline{s}}{M}\right)\nonumber
\end{equation}
for $L\mid n$. This completes the proof of the lemma.\\

We recall the notation and assumptions \eqref{Q} as well as the condition $(a, rs) = 1$. We apply the lemma with $(n, m) \leftarrow (an'', n')$, so that 
$$K = K_{d} = \Big(\frac{Q''}{(dQ_1, Q'')}, a^{\infty}\Big) \mid Q''$$
with the above notation. We conclude
$$\sqrt{an'n''}\rho_{f_d , 1/s}(\pm an'n'')  = \Lambda e\Big(\pm\frac{n'(an'' - K)\overline{s}}{rq}\Big)  \sqrt{K n' }\rho_{f_d , 1/s}(\pm K n' ) $$
where $\Lambda \ll (an'')^{\theta+\varepsilon}$ is independent of $n'$. We write $$a^{\ast}_{n'n''} = \delta_{(n', a) = 1}a_{n'n''}  e\Big(\pm\frac{n'(an'' - K)\overline{s}}{rq}\Big). $$
In this way 
the right hand side of \eqref{3-line} is bounded by
\begin{displaymath}
\begin{split}
(aNQ)^{\varepsilon}\sum_{q_0 \mid Q_1 \mid Q} & \sum_{d \mid \frac{Q}{Q_1}} \sum_{\substack{f \in \mathcal{B}^{\ast}(Q_1, Q, \chi)\\ |t_f| \leq T}}\frac{(1 + |t_f|)^{\pm \kappa} }{\cosh(\pi t_f)}   a^{2\theta} \\
& \times  \sum_{n'' \mid a^{\infty}} (n'')^{2\theta}    \Big|   \sum_{N/n'' < n' \leq 2N/n''} a^{\ast}_{n'n''}  \sqrt{K_dn'}\rho_{f_{d} , 1/s}(\pm K_dn')\Big|^2.
\end{split}
\end{displaymath}
As $f_{d}$ is $L^2$-normalized, we can complete it to an orthonormal basis of $\langle \{f_r : r \mid Q/Q_1\}\rangle$, restoring a full basis $\mathcal{B}(Q, \chi)$. 
 By positivity we can bound the previous display by 
$$ \ll (aNQ)^{\varepsilon} a^{2\theta} \sum_{n'' \mid a^{\infty}} (n'')^{2\theta}\sum_{K \mid Q''}  \sum_{\substack{f \in \mathcal{B}(Q, \chi)\\ |t_f| \leq T}}\frac{(1 + |t_f|)^{\pm \kappa} }{\cosh(\pi t_f)} \Big| \sum_{N/n'' < n \leq 2N/n''} a^{\ast}_{n'n''} \sqrt{Kn'} \rho_{f , 1/s}(\pm K  n')\Big|^2.$$
We summarize the discussion in the following lemma. 
\begin{lemma}\label{lem2-dr} Let $Q$ be as in \eqref{Q}, $(a, rs) = 1$. Then
$$ \mathcal{M}_{1/s}(T, Q, \chi, N, (a_n), a) \ll  (aNQ)^{\varepsilon} a^{2\theta} \sum_{\substack{n'' \mid a^{\infty}\\ n'' \leq 2N}} \sum_{K \mid Q''}  (n'')^{2\theta}
\mathcal{M}_{1/s}(T, Q, \chi, KN/n'', (a^{\ast}_n), 1) $$
where 
\begin{equation}\label{bn}
a^{\ast}_n = \delta_{(n, a) = 1} \delta_{K \mid n} a_{nn''/K} e\Big(\pm\frac{n(an'' - K)\overline{s}}{Krq}\Big). 
\end{equation}
which depends on $a, n'', K, rq, s$. The same inequalities hold for $\mathcal{H}_{1/s}$ and $\mathcal{E}_{1/s}$ in place of $\mathcal{M}_{1/s}$ except that the factors $a^{2\theta}$ and $(n'')^{2\theta}$ can be removed.\end{lemma}

\textbf{Proof.} The case of $\mathcal{H}_{1/s}$ is completely analogous. Indeed all the arguments carry through after replacing $f(g.z)$ by $ \det(g)^{k/2} j(g,z)^{-k}f(gz)$ wherever necessary. In particular, as mentioned in \cite{BM} just below Lemma~2, if we define $f_d(z) = d^{k/2}f(dz)$, then \eqref{0new} remains valid and $\Vert f_d\Vert = \Vert f\Vert$ with respect to the Petersson norm. Here we can even take $\theta=0$, since the Ramanujan conjecture is a theorem due to Deligne.

We now turn towards $\mathcal{E}_{1/s}$.  Here the key difficulty is to choose a basis for the continuous spectrum such that Lemma~\ref{lem1-dr} remains true for the Fourier coefficients of the Eisenstein series. All the necessary ingredients are provided in \cite{Y}.

Given two Dirichlet characters $\chi_1$ and $\chi_2$, we associate the Eisenstein series $E_{\chi_1,\chi_2}(z,s)$ as in \cite[(3.3)]{Y}. These Eisenstein series are eigenfunctions of all Hecke operators (see \cite[(4.14)]{Y}) and their Fourier expansion has been computed in \cite[Proposition~4.1]{Y}. Let $q_1$ be the conductor of $\chi_1$ and $q_2$ the conductor of $\chi_2$. The space of Eisenstein series of level $Q$ and nebentypus $\chi$ is spanned by
\begin{equation}
	\{ E_{\chi_1,\chi_2}(dz,s)\colon \chi_1\chi_2^{-1} = \chi,\, q_1q_2d\mid Q\}.\nonumber
\end{equation}
As discussed in \cite[Section~8.5]{Y} one can extract an orthogonal basis as above. To be more precise, an orthogonal basis for the space of Eisenstein series of level $Q$ and nebentypus $\chi$ is given by
\begin{equation}
	\bigsqcup_{Q_1\mid Q}\bigsqcup_{\substack{\chi_i \text{ mod }q_i\text{ primitive},\\ q_1q_2 = Q_1, \\ \chi_1\chi_2^{-1}=\chi}} \left\{ E_{\chi_1,\chi_2}^{(g)}(z,s)\colon g\mid \frac{Q}{Q_1}\right\}.\nonumber
\end{equation}
Furthermore \eqref{0new} holds with $\theta=|\Re(s)-1/2|$. This enables us to reproduce the argument above line by line as long as Lemma~\ref{lem1-dr} remains true. Looking into the proof shows that it  only relies on elementary matrix manipulations and the fact that newforms are Atkin-Lehner pseudo-eigenfunctions. Fortunately the latter is also true for the Eisenstein series $E_{\chi_1,\chi_2}(z,s)$, as shown in \cite[Section~9.1]{Y}. We conclude that the continuous contribution can be handled exactly as the Maa\ss\  form contribution detailed above. \\

 We conclude by \cite[Proposition 4.7]{Dr} that for $\mathfrak{a} = 1/s$ each of the three spectral quantities \eqref{spectral} is
\begin{equation}\label{47}
\ll (QNa)^{\varepsilon} a^{2\theta} \sum_{\substack{n'' \mid a^{\infty}\\ n'' \leq 2N}}  \sum_{K \mid Q''}  (n'')^{2\theta}\Big(T^2 + q_0^{1/2} \frac{1}{Q} \frac{KN}{n''}\Big) \| \textbf{a}^{\ast} \|_2^2. 
\end{equation}
If $a= 1$ (which will be the case in our application for the term corresponding to the cusp $\infty$), we can simply use  the standard large sieve inequality as stated in \cite[Proposition 4.7]{Dr}. 

To treat the exceptional spectrum we define
$$E_{q, \mathfrak{a}}(Y, (a_n), N, a) =  \sum_{\substack{f \in \mathcal{B}(q, \chi)\\ t_f \in i\Bbb{R}}} Y^{2|t_f|} \Big| \sum_{N < n \leq 2N} a_n \sqrt{an} \rho_{f , \mathfrak{a}}(an)\Big|^2$$
for $Y \geq 1$. By the same argument as in Lemma \ref{lem2-dr} we have
$$E_{q, 1/s}(Y, (a_n), N, a) \ll  (QNa)^{\varepsilon}a^{2\theta}\sum_{\substack{n'' \mid a^{\infty}\\ n'' \leq 2N}}  \sum_{K \mid Q''}  (n'')^{2\theta} E_{q, 1/s}(Y, (a^{\ast}_n),KN/n'', 1)$$
with $a^{\ast}_n$ as in \eqref{bn}. By \cite[Lemma 4.8]{Dr} we conclude
\begin{equation}\label{48}
\begin{split}
&E_{q, 1/s}(Y, (a_n), N, a)\\
& \ll  (QNa)^{\varepsilon}a^{2\theta}\sum_{\substack{n'' \mid a^{\infty}\\ n'' \leq 2N}}  \sum_{K \mid Q''}  (n'')^{2\theta} \Big(1 + \Big(\frac{1}{Q} \frac{KN}{n''}Y\Big)^{1/2}\Big) \Big(1 + \Big(\frac{q_0}{Q} \frac{KN}{n''}\Big)^{1/2}\Big)  \| \textbf{a}^{\ast} \|_2^2. 
\end{split}
\end{equation}
We will use \cite[Lemma 4.10]{Dr} without any modifications for the cusp $\infty$. 

\subsection{Bounds for sums of generalized Kloosterman sums} The scene has now been prepared to complete the proof of Theorem \ref{thm-dr}. We start with the following modification of \cite[Proposition 4.12]{Dr}.
\begin{lemma}\label{lem3-dr} Let $M, N, R, S \geq 1$, $X > 0$, $a, q\in \Bbb{N}$, $\chi$ a Dirichlet character modulo $q$, $\Phi$  a smooth function supported on $[X, 2X]$ such that $\| \Phi^{(j)}\|_{\infty} \ll_j X^{-j}$, let $(a_m)$, $(b_{n, r, s})$ be sequences of complex numbers supported on $M < m \leq 2M$, $N < n \leq 2N$, $R < r \leq 2R$, $S < s \leq 2S$. Assume that $a_m$ is the characteristic function of an interval and $b_{n, r, s}$ vanishes unless $(rs, a) = 1$. Then
\begin{displaymath}
\begin{split}
&\sum_{\substack{m, n, r, s\\ (s, rq) = 1}} a_m b_{n, r, s} \sum_{c \in \mathcal{C}(\infty, 1/s)} \frac{1}{c} \Phi\Big(\frac{4\pi \sqrt{amn}}{c}\Big) S_{\infty, 1/s}(m, \pm an, c) \ll (aq(X+X^{-1})RSMN)^{\varepsilon} \big(L_1 + L_2\big)
\end{split}
\end{displaymath}
where
\begin{displaymath}
\begin{split}
& L^2_1 =  a^{2\theta} \sum_{\substack{n'' \mid a^{\infty}\\ n'' \leq 2N}}  \sum_{K \mid Q''}  (n'')^{2\theta} \Big(1 + X + \sqrt{\frac{N'}{RS}}\Big)^2\Big(1 + X + \sqrt{\frac{M}{RS}}\Big)^2 RS \frac{   M  }{1+X^2} \|  \tilde{\textbf{b}}(n'')  \|^2_2 \\
& L_2^2 \ll  a^{2\theta} \sum_{\substack{n'' \mid a^{\infty}\\ n'' \leq 2N}}  \sum_{K \mid Q''}  (n'')^{2\theta}\Big(1 + \sqrt{\frac{N'}{RS}}\Big)^2  \Big(\frac{(N'M)^{1/4}(1 + X^{-1/2})}{(N' + RS)^{1/4}} \Big)^2\frac{M }{1+X^2}\| \tilde{\textbf{b}}(n'')  \|_{2}^2
\end{split}
\end{displaymath}
with 
\begin{equation}\label{btilde}
N' = \frac{KN}{n''}, \quad \text{and} \quad \tilde{b}_{n, r, s}(n'') =   b_{nn'' , r, s} ,
\end{equation}
and the Kloosterman sums is defined with respect to the group $\Gamma_0(qrs)$ as in \cite[Section 4]{Dr}. 
\end{lemma}

\textbf{Proof.} The bound $L_1$ for the regular spectrum follows without difficulty from the large sieve inequality, both in its standard version \cite[Proposition 4.7]{Dr} and in our modified version \eqref{47}. The bound for $L_2$ follows as  in \cite{Dr}, based on \cite[Section 9.1]{DI}. As the argument are not completely straightforward and \cite{DI} contains several typos, we present some details. By the Kuznetsov formula and the Cauchy-Schwarz inequality we need to bound $(1+X)^{-1}( \mathcal{M}_1\mathcal{M}_2)^{1/2}$ where
\begin{displaymath}
\begin{split}
& \mathcal{M}_1 = \sum_{(s, rq) = 1} \sum_{\substack{f \in \mathcal{B}(qrs, \chi)\\ t_f \in i \Bbb{R}}} \Big( \frac{1 + X^{-2}}{1 + X_0}\Big)^{2|t_f|} \Big|\sum_m a_m m^{1/2} \rho_{f, \infty}(m)\Big|^2,\\
&  \mathcal{M}_2 = \sum_{(s, rq) = 1} \sum_{\substack{f \in \mathcal{B}(qrs, \chi)\\ t_f \in i \Bbb{R}}} (1 + X_0)^{2 |t_f|} \Big| \sum_n b_{n, r, s} n^{1/2} \rho_{f, 1/s}(an)\Big|^2.
\end{split}
\end{displaymath}
and $X_0>0$ is a free parameter at our disposal. We choose $X_0 = qRS/N'$. By \eqref{48} we obtain
$$ \mathcal{M}_2 \ll (QNa)^{\varepsilon}a^{2\theta} \sum_{\substack{n'' \mid a^{\infty}\\ n'' \leq 2N}}  \sum_{K \mid Q''}  (n'')^{2\theta}  \Big(1 + \frac{N'}{RS}\Big)   \| \textbf{b}^{\ast}(n'', K) \|_{2}^2$$
where
$$b^{\ast}_{n, r, s}(n'', K) = \delta_{(n, a) = 1} \delta_{K \mid n} b_{nn''/K, r, s} e\Big(\pm\frac{n(an'' - K)\overline{s}}{Krq}\Big),$$
so that  $ \| \textbf{b}^{\ast}(n'', K) \|_{2} \leq \| \tilde{\textbf{b}}(n'') \|_2$. 
We treat $\mathcal{M}_1$ by   \cite[Lemma 4.10]{Dr} getting
$$\mathcal{M}_1 \ll (qRSM)^{\varepsilon}\Big(RS + M +  M^{1/2}\Big(1 +  \frac{1 + X^{-2}}{1 + qRS/N'} \Big)^{1/2}\Big)M.$$
We can drop the first two terms in the parentheses as their contribution is already covered by $L_1$, and the result follows.\\ 

This leads to the following version of \cite[Proposition 4.13]{Dr} whose proof carries over verbatim. 
\begin{lemma}\label{lem4-dr} Let $M, N, R, S, C \geq 1$, $a, q \in \Bbb{N}$,   $\chi$ a Dirichlet character modulo $q$, $g$ a smooth function supported on $[C, 2C] \times [M, 2M] \times \Bbb{R}_{>0}^3$ such that 
$$\frac{\partial^{\nu_0 + \ldots + \nu_4}}{\partial c^{\nu_0} \partial m^{\nu_1} \partial n^{\nu_2} \partial r^{\nu_3} \partial s^{\nu_4} } g(c, m, n, r, s) \ll _{{\bm \nu}} C^{-\nu_0} M^{-\nu_1} N^{-\nu_2} R^{-\nu_3} S^{-\nu_4} $$
for ${\bm \nu} \in \Bbb{N}_0^5$. Let $(b_{n, r, s})$ be sequences of complex numbers supported on   $N < n \leq 2N$, $R < r \leq 2R$, $S < s \leq 2S$. Assume that   $b_{n, r, s}$ vanishes unless $(rs, a) = 1$.  Let $t \in [0, 1]$. Then
\begin{displaymath}
\begin{split}
& \sum_{\substack{c, m, n, r, s\\ (sc, rq) = 1}} b_{n, r, s} \bar{\chi}(c) g(c, m, n, r, s)e(mt) S(an\bar{r}, \mp m\bar{q}, sc)  \ll (CRSMNqa)^{\varepsilon} q^{3/2}(K_1 + K_2) \sqrt{M }  
\end{split}
\end{displaymath}
where
\begin{displaymath}
\begin{split}
&K_1^2 = a^{2\theta}\sum_{\substack{n'' \mid a^{\infty}\\ n'' \leq 2N}}  \sum_{K \mid Q''}  (n'')^{2\theta} RS \frac{(C^2 S^2 R + aMN + C^2 S N')(C^2 S^2 R + aMN + C^2 S M)}{C^2 S^2 R + aMN} \| \tilde{ \textbf{b}}(n'')\|_2^2,\\
& K_2^2 =a^{2\theta}\sum_{\substack{n'' \mid a^{\infty}\\ n'' \leq 2N}}  \sum_{K \mid Q''}  (n'')^{2\theta}  C^3 S^2 \sqrt{R(RS + N')} \Big(\frac{K}{an''}\Big)^{1/2} \|  \tilde{ \textbf{b}}(n'') \|_2^2
\end{split}
\end{displaymath}
and $N'$ and $\tilde{ \textbf{b}}(n'')$ have the same meaning as in \eqref{btilde}. 
\end{lemma}


To complete the proof of Theorem \ref{thm-dr}, we follow \cite[Section 4.3.3]{Dr}. Note that there is a small error in the display after in \cite[(4.36)]{Dr} for $\mathcal{A}_0$, which is already present in \cite[(9.11)]{DI}. This was recently corrected in \cite{BFI3}, namely the gcd $(n, cs)$ cannot be estimated trivially and the factor $S^{-1/2}$ should be removed. The rest of the proof follows as in \cite[p.\ 707]{Dr} by using the following inequalities with the notation in Lemma \ref{lem4-dr}. If $D\geq 1$ is another parameter and  $M \ll SCq/D$, we have
\begin{displaymath}
\begin{split}
\frac{D^2M}{(SC)^2} K_1^2 &\ll a^{2\theta}\sum_{\substack{n'' \mid a^{\infty}\\ n'' \leq 2N}}  \sum_{K \mid Q''}  (n'')^{2\theta}  \frac{D^2M}{(SC)^2} RS\Big( C^2 S^2 R +aMN + \frac{ C^2MN'}{R}+  C^2S(N'+M)\Big)  \|\tilde{ \textbf{b}}(n'')\|_2^2\\
& \ll q^2 (Q'')^{1 + \varepsilon} a^{2\theta}\sum_{\substack{n'' \mid a^{\infty}\\ n'' \leq 2N}}   (n'')^{2\theta} \Big(CS\Big(\frac{N}{n''} + RS\Big)(C + RD)  + a SNR\Big)  \|\tilde{ \textbf{b}}(n'') \|_2^2. 
\end{split}
\end{displaymath}
and
\begin{displaymath}
\begin{split}
\frac{D^2M}{(SC)^2} K_2^2 &  \ll  q (Q'')^{1 + \varepsilon}  C^2 D S \sqrt{R(RS + N)}  \| \tilde{ \textbf{b}}(n'') \|_2^2. 
\end{split}
\end{displaymath}
These bounds lead directly to  the bounds in Theorem \ref{thm-dr} by using  $\| \tilde{ \textbf{b}}(n'') \|_2 \leq \| \textbf{b}\|_2$ and $Q'' \leq q$.


\section{Proof of  Theorem \ref{PrimeAp}}\label{sec4}

For $R \geq 1$, $n \in \Bbb{Z} \setminus \{0\}$, $q\in \Bbb{N} $ denote 
$$u_R(n;q)=\mathds{1}_{n\equiv 1 \, (\text{mod } q)}-\frac{1}{\phi(q)}\sum_{\substack{\chi\, (\text{mod  }q)\\ \operatorname{cond}(\chi)\leq R}}\chi(n)=\frac{1}{\phi(q)}\sum_{\substack{\chi\, (\text{mod } q)\\ \operatorname{cond}(\chi)>R}}\chi(n).$$
Note that $u_R(n;q)=0$ if $(q, n)>1$ or $q\leq R$. We have the trivial bound 
$$u_R(n;q)\ll \mathds{1}_{n\equiv 1\, (\text{mod } q)}+ \frac{\tau(q)}{\phi(q)}R.$$
The following is an analogue of \cite[Theorem 5.1]{Dr} with a longer range for $a_1$ and an additional congruence condition on $q$.

\begin{prop}\label{TypeII}
	Let $M, N, Q, R\geq 1$, $a_1, a_2 \in \mathbb{Z}\backslash \{0\}$, write  $x=MN$. Let $c\in \Bbb{N}$, $c_0 \in \Bbb{Z}$ with $(c_0, c) = 1$. Let $\alpha_m$ and $\beta_n$ be two sequences supported in $m\in (M, 2M]$ and $n\in (N, 2N]$ such that 
	$\alpha_m\leq \tau(m)^A, \beta_n\leq \tau(n)^A$ 
	for some $A\geq 1.$ 
	Let $\eta > 0$ be any sufficiently small number. Then there exist $\delta = \delta(\eta) > 0$ and $D = D(\eta, A)$ with the following property.
		If	
  \begin{equation*}
    x^\eta\leq N \leq x^{1/4+\eta},  \quad 
     Q\leq x^{1/2+\delta}, \quad 
         c, R, |a_2|\leq x^\delta, \quad 
     |a_1|\leq x^{1+\delta} 
     \end{equation*}
	then	
$$	\sum_{\substack{Q\leq q\leq 2Q\\ (q, a_1a_2)=1\\ q\equiv c_0 \,(\text{{\rm mod }}  c)}}\sum_{\substack{m, n\\(n, a_2)=1}}\alpha_m\beta_{n}u_R(mn\overline{a_1}a_2;q)\ll_{ \eta, A} cx(\log x)^{D}R^{-1}.$$
\end{prop}
If $Q\leq R$, we have $u_R(n;q)=0$. For $R\leq Q\leq x^{1/2}/R$, we can apply the Bombieri-Vinogradov theorem as in \cite[Lemma 5.2]{Dr}. Thus it remains to deal with the case 
  \begin{equation}\label{bounds}
    x^\eta\leq N \leq x^{1/4+\eta},  \quad 
     x^{1/2-\delta} \leq Q\leq x^{1/2+\delta}, \quad 
         c, R, |a_2|\leq x^\delta, \quad 
     |a_1|\leq x^{1+\delta}. 
     \end{equation}
Verbatim as in \cite[pp.\ 709-710]{Dr}, we may assume that $\beta_n$ is supported on squarefree integers (by estimating large square factors trivially and moving small square factors into $a_2$) and we can also smooth out the $q$-sum. In this way, Proposition \ref{TypeII} follows easily from the following statement: 
\begin{prop}\label{TypeIIProp}
Let $x, M, N, Q, R, \eta, (\alpha_m), (\beta_n)$ be as in Proposition \ref{TypeII}. 
Assume in addition that $\beta_n$ is supported on squarefree  integers $n$. Then for any smooth function $\gamma:\mathbb{R}_+\rightarrow [0,1]$ with 
	$$\mathds{1}_{Q\leq q\leq 2Q}\leq \gamma(q)\leq \mathds{1}_{Q/2\leq q\leq 3Q/2}$$
	and $\|\gamma^{(j)}\|_\infty\ll_j Q^{-j+B\delta j}$ for some $B\geq 0$ and for all  $j\geq 0$  we have 
$$\mathcal{D} := \sum_{\substack{ (q, a_1a_2) = 1\\ q\equiv c_0 \, (\text{{\rm mod }} c)}}\gamma(q)\sum_{\substack{m, n\\ (n, a_2)=1}}\alpha_m \beta_nu_R(mn \overline{a_1}a_2;q)\ll cx(\log x)^{O(1)}R^{-1}. $$
The implicit constants depend on $\eta, A$ and $B$ at most.
\end{prop}
\textbf{Proof.} 	The proof follows along the   lines of the proof of  \cite[Proposition 5.3]{Dr}, to some extent verbatim, the only differences being the longer range of $a_1$ and the additional congruence condition on $c$. We highlight the points where some extra input is needed. In the following, the value of $\delta$ may change from line to line. 

We start by observing that  the terms from $m,n \mid a_1$ contribute to $\mathcal{D}$ at most 
	$$\sum_{q}\gamma(q)\sum_{m\mid a_1, n\mid a_1}|\alpha_m\beta_n|\Big (1+ \frac{R\tau(q)}{\phi(q)}\Big)\ll x^\epsilon (Q R)\ll x^{1/2+\delta+\epsilon}, $$
which is acceptable. 	So we can restrict $m, n$ such that $mna_2-a_1\not=0$. This is important in the display after \cite[(5.23)]{Dr}. By triangle inequality, we write 
	\begin{align*}
	\mathcal D\ll \sum_{b\, (\text{mod } c)}\sum_{m}\Big( |\alpha_m|\Big|\sum_{q}\sum_{n\equiv b \, (\text{mod } c)} (\ldots) \Big|\Big).
	\end{align*}
	Thus, without loss of generality we can assume that $\beta_n=0$ unless $n\equiv b_0$ (mod $c$) for some $1\leq b_0\leq c$ and recover the bound for $\mathcal{D}$ for general $\beta_n$  by a factor of $c$. After an application of the Cauchy-Schwarz inequality, we arrive at quantities $\mathcal{S}_i, i=1,2,3$ as in \cite[(5.13)]{Dr} with the difference being the additional congruence condition on $q$, which can be carried through in the evaluations of $\mathcal{S}_2$ and $\mathcal{S}_3$ in \cite[Section 5.3.1/5.3.2]{Dr}. We can also reduce $\mathcal{S}_1$ to $\mathcal{S}_1(q_0, n_0)$ as in \cite[Section 5.4]{Dr}.

	The interesting part happens in \cite[Section 5.5]{Dr}. Here we need to evaluate
	$$\mathcal{S}_1(q_0, n_0)=\sum_{\substack{  (q_1, q_2)=1\\ (q_1q_2, a_1a_2)=1\\ q_0q_j \equiv c_0 \, (\text{mod } c)}}\gamma(q_0q_1)\gamma(q_0q_2)\sum_{\substack{ (n_1, n_2) = 1\\ (n_0n_j, q_0q_ja_2)=1\\ n_1\equiv n_2 \,(\text{mod } q_0)}}\beta_{n_0n_1}\overline{\beta_{n_0n_2}}\sum_{\substack{m\equiv a_1\overline{a_2n_0n_j} \, (\text{mod } q_0q_j) }}\alpha(m)$$
where $q_0, n_0 \leq x^{\delta}$, $(q_0, a_1a_2, c) = (n_0, a_2) = 1$ and $\alpha(m)$ is a smooth function with support in $m \asymp M$. 
	We use M\"obius inversion to detect the condition $(q_1q_2, a_1)=1$ getting
	\begin{align*}
	&\mathcal{S}_1(q_0, n_0)\\
	&=\sum_{\substack{\delta_1 \delta_2\mid a_1\\ (\delta_1, \delta_2)=1}} \mu(\delta_1\delta_2) \sum_{\substack{ (\delta_1q_1, \delta_2q_2)=1\\ (\delta_1\delta_2q_1q_2, a_2)=1\\q_0\delta_jq_j\equiv c_0 \, (\text{mod } c)}}\hspace{-0.2cm} \gamma(q_0\delta_1q_1)\gamma(q_0\delta_2q_2)\hspace{-0.2cm}\sum_{\substack{n_1, n_2\\ (n_0n_j, q_0\delta_jq_ja_2)=1}}\hspace{-0.2cm}\beta_{n_0n_1}\overline{\beta_{n_0n_2}}\hspace{-0.3cm}\sum_{\substack{m\equiv a_1\overline{a_2n_0n_j} \, (\text{mod }q_0\delta_jq_j)}}\hspace{-0.2cm}\alpha(m)\\
	&:=\sum_{\substack{\delta_1, \delta_2\mid a_1\\ (\delta_1, \delta_2)=1}}\mu(\delta_1\delta_2)\mathcal{S}_1(q_0, n_0, \delta_1, \delta_2),
	\end{align*}
say. 	Note that trivially 
	\begin{align*}
\mathcal{S}_1(q_0, n_0, \delta_1, \delta_2) 
	& \ll x^\epsilon \sum_{\substack{(\delta_1\delta_2q_1 q_2, a_2)=1\\q_1\asymp Q/q_0\delta_1 }} \sum_{ \substack{n_j \asymp N/n_0\\ n_1\equiv n_2\,(\text{mod } q_0)} } \sum_{\substack{ a_2mn_0n_1\equiv a_1 \, (\text{mod }q_0\delta_jq_j)\\ a_2mn_0n_2\not=a_1}}1\ll  x^\epsilon N\Big(\frac{MN}{\max(\delta_1, \delta_2)}+1\Big).
	\end{align*}
		Therefore, the contribution of the terms with $\max(\delta_1, \delta_2)>x^{\delta}$ can be bounded by $x^\epsilon MN^2 x^{-\delta}$.  
The same analysis works for the expected main term, and   it is enough to show for $$\delta_1, \delta_2\leq x^\delta, \quad (q_0, a_1a_2c) = (n_0, a_2) = (\delta_1\delta_2, a_2c) = 1$$ that
	\begin{align*}
	\mathcal{S}_1(q_0, n_0, \delta_1, \delta_2)=\hat{\alpha}(0)X_1(q_0, n_0, \delta_1, \delta_2)+O(MN^2x^{-\delta}).
	\end{align*}
	where 
$$X_1(q_0, n_0, \delta_1, \delta_2) = 	\sum_{\substack{(q_1q_2, a_2)=1\\ (\delta_1q_1,\delta_2q_2)=1\\q_0\delta_j q_j \equiv c_0  \, (\text{mod } c)}}\frac{\gamma(q_0\delta_1 q_1)\gamma(q_0\delta_2q_2)}{q_0\delta_1q_1\delta_2q_2}\sum_{\substack{(n_1, n_2)=1\\(n_0n_j, q_0\delta_jq_ja_2)=1\\  n_1\equiv n_2\,(\text{mod } q_0)}}\beta_{n_0n_1}\overline{\beta_{n_0n_2}}.$$
	After Poisson summation, we need to bound the non-zero frequencies
	$$\mathcal{R}_1=\sum_{\substack{(q_1 q_2, a_2) = 1\\ (\delta_1q_1, \delta_2q_2)=1\\ q_0\delta_jq_j 
	\equiv c_0\, (\text{mod } c)}}\gamma(q_0\delta_1q_1)\gamma(q_0\delta_2q_2)\sum_{\substack{ (n_1, n_2)=1\\ (n_0n_j, q_0\delta_jq_ja_2)=1\\ n_1\equiv n_2\, (\text{mod } q_0)}}\beta_{n_0n_1}\overline{\beta_{n_0n_2}}\sum_{0<|h|\leq H}\frac{1}{W}\hat{\alpha}\Big( \frac{h}{W}\Big)e\Big(\frac{h\mu}{W}\Big)$$
where
$$	
 W=q_0\delta_1q_1\delta_2q_2, \quad  H=W^{1+\epsilon}M^{-1}$$  	and $\mu$ satisfies $\mu \equiv a_1\overline{a_2n_0n_j}\, (\text{mod } {q_0\delta_jq_j})$ for $j=1, 2$. 
Since $(n_0, n_1)=1$, we have
\begin{align*}
\frac{\mu}{q_0\delta_1q_1\delta_2q_2}\equiv \frac{a_1}{q_0\delta_1q_1\delta_2q_2a_2n_0n_1}+a_1\frac{n_1-n_2}{q_0}\frac{\overline{\delta_1q_1a_2n_0n_2}}{n_1\delta_2q_2}-a_1\frac{\overline{q_0\delta_1q_1\delta_2q_2n_1}}{a_2n_0}\pmod {1}.
\end{align*}
From
\begin{align*}
\frac{a_1}{q_0\delta_1q_1\delta_2q_2a_2n_0n_1}\ll\frac{a_1}{x} \frac{x}{WN}\ll x^\epsilon\frac{a_1}{x} H^{-1}
\end{align*}
we see that first term in the exponential is  essentially ``flat'' and can be removed by partial summation, so that 
 $$\mathcal{R}_1\ll x^\epsilon \left(\frac{a_1}{x}+1\right)\sup_{\substack{Q\leq Q'\leq 2Q\\ N\leq N'\leq 2N\\ 1\leq H'\leq H}}| \mathcal{R}_1'|,$$ where
\begin{align*}
&\mathcal{R}_1'= \sum_{\substack{(q_1 q_2, a_2) = 1\\ (\delta_1q_1, \delta_2q_2)=1\\ q_0\delta_jq_j 
\equiv c_0 \, (\text{mod } c)\\ \delta_jq_jq_0\leq Q'}}{\gamma}(q_0\delta_1q_1) {\gamma}(q_0\delta_2q_2) \sum_{\substack{ (n_1, n_2)=1\\ (n_0n_j, q_0\delta_jq_ja_2)=1\\ n_1\equiv n_2 \, (\text{mod } q_0)\\ n_0n_j\leq N'}}{\beta}_{n_0n_1}\overline{{\beta}_{n_0n_2}}\\& \times\sum_{0<|h|<H'}\frac{1}{W}\hat{\alpha}\Big(\frac{h}{W}\Big)e\Big(a_1\frac{n_1-n_2}{q_0}\frac{\overline{\delta_1q_1a_2n_0n_2}}{n_1\delta_2q_2}-a_1\frac{\overline{q_0\delta_1q_1\delta_2q_2n_1}}{a_2n_0}\Big).
\end{align*}
We incorporate the condition $\delta_jq_1q_0\leq Q'$ in the smooth function $\gamma$ without changing the notation. Similarly,  we put condition $n_0n_j\leq N'$ in $\beta_n$ without changing the notation. We then obtain
\begin{align*}
\mathcal{R}_1'\ll_\epsilon x^\epsilon (n_0|a_2|)^2\frac{Mq_0}{Q^2}\sup _{\xi \asymp  Mq_0/Q^2}\max_{\substack{\lambda_1, \lambda_2 \, (\text{mod } n_0a_2)\\  (\lambda_1\lambda_2, n_0a_2) = 1\\\lambda_j\equiv c_0 \overline{q_0} \, (\text{mod } (n_0a_2, c))}}|\mathcal{R}_1''|,
\end{align*} 
where 
\begin{align*}
\mathcal{R}_1''=&\sum_{\substack{  (\delta_1q_1, \delta_2q_2)=1\\ q_j \equiv \lambda_j\overline{\delta_j}\, (\text{mod } n_0a_2)\\ q_j \equiv c_0\overline{\delta_jq_0} \, (\text{mod } c)}}\gamma(q_0\delta_1q_2)\gamma(q_0\delta_2q_2)\sum_{\substack{n_1, n_2\\ (n_1, n_2)=1\\ (n_0n_j, q_0q_j\delta_ja_2)=1\\ n_1\equiv n_2\, (\text{mod } q_0)}}\beta_{n_0n_1}\overline{\beta_{n_0n_2}}\\
& \times \sum_{0<|h|< H'}\alpha ( \xi q_0\delta_1\delta_2q_1q_2)e\Big(-\xi h-a_1h\frac{\overline{q_0\lambda_1\lambda_2n_1}}{a_2n_0}\Big)e\Big(a_1h\frac{n_1-n_2}{q_0}\frac{\overline{q_1a_2n_0n_2\delta_1}}{n_1\delta_2q_2}\Big). 
\end{align*}

Now we can apply Theorem \ref{thm-dr} with 
\begin{displaymath}
\begin{split}
&(\boldsymbol c, \boldsymbol d, \boldsymbol a, \boldsymbol n, \boldsymbol r, \boldsymbol s,\boldsymbol q)\leftarrow (q_2, q_1,a_1,h\frac{n_1-n_2}{q_0}, a_2n_0n_2\delta_1, n_0\delta_2, [n_0a_2,c]),\\
&(\boldsymbol C,\boldsymbol D,\boldsymbol N,\boldsymbol R, \boldsymbol S)\leftarrow \Big(\frac{Q}{q_0\delta_2}, \frac{Q}{q_0\delta_1}, \frac{HN}{q_0n_0}, a_2\delta_1N, \frac{N\delta_2}{n_0}\Big),\\
&g(\textbf{c}, \textbf{d}, \boldsymbol{n}, \boldsymbol{r}, \boldsymbol{s})= \gamma(q_0\delta_1 \textbf{c})\gamma(q_0\delta_2\textbf{d})\alpha(\xi q_0\delta_1\delta_2\textbf{c} \textbf{d}),
\end{split}
\end{displaymath}
and 
\begin{align*}
\tilde{b}_{\boldsymbol{n},\boldsymbol{r},\boldsymbol{s}}{(n'')}=\operatornamewithlimits{\sum\sum}_{\substack{(n_1,n_2) = 1\\ (n_0n_j, q_0\delta_jq_ja_2)=1\\ n_1\equiv n_2\, (\text{mod } q_0)\\ \boldsymbol{r}=a_2n_0n_2\delta_1\\ \boldsymbol{s}=n_1\delta_2}}\beta_{n_0n_1}\overline{\beta_{n_0n_2}}\sum_{\substack{0<|h|< H'\\q_0n''\boldsymbol{n}=h(n_1-n_2) }}\alpha(\xi q_0\delta_1\delta_2q_1q_2)e\Big(-\xi h-a_1h \frac{\overline{q_0\lambda_1\lambda_2n_1}}{a_2n_0}\Big).
\end{align*}
We now verify the condition $(\boldsymbol{q}, \boldsymbol{cs})=1$. We have $(n_0a_2, q_2n_1\delta_2)=1$ from $(n_0n_j, q_0q_j\delta_ja_2)= (n_0, n_1) = (q_j\delta_j, n_0a_2)=1$ and we have $(c, q_2\delta_2)=1$ since $(c_0, c)=1$. By the assumption that $\beta_n$ vanishes unless $n\equiv b_0 \, (\text{mod }c)$, we have $n_0n_1\equiv n_0n_2\, (\text{mod }c)$ and thus we have $(b_0, c)\mid n_0$. Since $\beta_n$ is supported on squarefree integers, we have $(n_1, c)=1$.  Since $(n_1, n_2)=1$ and $N\geq x^\eta$, we have $n_1\not=n_2$ and 
\begin{align*}
\|\tilde{b}_{\boldsymbol{n},\boldsymbol{r},\boldsymbol{s}}(n'')\|_2^2&\ll \sum_{n_1}\sum_{n_2}|\beta_{n_0n_1}\beta_{n_0n_2}|\sum_{\substack{\frac{n''}{(n'', n_1-n_2)}\mid h\leq H'} }1\ll\frac{H'}{n''}\sum_{n_1, n_2}(n_1-n_2, n'')\|\beta_n\|_\infty^2\\
&\ll \frac{H'}{n''}\sum_{d\mid n''}d N\Big(\frac{N}{d}+1\Big)\|\beta_n\|_\infty^2\ll\Big( \frac{H'}{n''}N^2+ H'N \Big)\tau(n'')\|\beta_n\|_\infty^2, 
\end{align*}
thus
\begin{align*}
\sum_{\substack{n''\mid a_1^\infty \\ n''\leq 2N}}(n'')^{2\theta}\|\tilde{b}_{\boldsymbol{n},\boldsymbol{r},\boldsymbol{s}}(n'')\|_2^2\ll H'N^2\tau(a_1) \|\beta_n\|_\infty^2\ll x^\epsilon H'N^2.
\end{align*}
We conclude
 $R_1''\ll x^{O(\delta)}\mathcal{A}^{1/2} \mathcal{B}^{1/2}$, where 
\begin{align*}
\mathcal{A}&\ll HN^2, \quad \mathcal{B}\ll a_1^{2\theta}( Q^2 N^3(N+ H)+a_1HN^3)+(Q^3 N^2\sqrt{N+H})+Q^2HN^2.
\end{align*}
Since  $x = MN$ and $Q \leq x^{1/2 + \delta}$, we have
 $$H=\frac{W^{1+\epsilon}}{M}\ll x^{O(\delta)} \frac{Q^2}{M} \ll x^{O(\delta)}N$$
and thus $\mathcal{B} \ll (a_1^{2\theta}Q^2N^4 + Q^3N^{5/2}) x^{O(\delta)}$  
 for  $N, Q, a_1$ satisfying \eqref{bounds}.  
Therefore,   
\begin{align*}
\mathcal{R}_1&\ll \left(\frac{a_1}{x}+1\right)\frac{M}{Q^2}x^{O(\delta)}\mathcal{A}^{1/2}\mathcal{B}^{1/2}\ll x^{O(\delta)} (a_1^{\theta} N^{3/2} Q^{-1}  +  N^{3/4} Q^{-1/2}) MN^2\\
& \ll (x^{\theta-\frac{1}{8}+3\eta/2+O(\delta)} +x ^{3/4\eta -1/16+O(\delta)})M N^2.
\end{align*}
Provided that $\theta < 1/8$ and $\eta$ is sufficiently small, this satisfies the required bound $O(MN^2x^{-\delta})$
by taking $\delta$ small enough depending on $\eta$.
The rest follows the same way as in \cite[Section 5.6]{Dr} as we can drop the congruence conditions on $q_i$ as an upper bound for $X_1-X_3$. This completes the proof of Proposition \ref{TypeIIProp} and hence  of Proposition \ref{TypeII}. \\

With Proposition \ref{TypeII} at hand, we can now proceed to the \textbf{proof of Theorem \ref{PrimeAp}}. Let $R\leq x^{\delta}$ be a parameter. We will choose $R = (\log x)^A$ unconditionally and $R = x^{\delta}$ on GRH. Note that the conditions imply $(d, q)= 1$. We write
\begin{align*}
&\sum_{\substack{q\leq Q\\ (q, a_1a_2)=1\\ q\equiv c_0\, (\text{mod }c)}}\sum_{\substack{n\leq x\\ n\equiv a_1 \overline{a_2}\, (\text{mod } q)\\ n\equiv d_0 \, (\text{mod } d)}}\Lambda(n)=\sum_{\substack{q\leq Q\\ (q, a_1a_2)=1\\ q\equiv c_0\, (\text{mod }c)}}\frac{1}{\phi(q)\phi(d)}\sum_{\substack{\psi \, (\text{mod }d)\\ \chi\, (\text{mod } q)}}\sum_{n\leq x}\Lambda(n)\chi(n\overline{a_1}a_2)\psi(n\bar{d}_0)\\
&=\sum_{\substack{q\leq Q\\ (q, a_1a_2)=1\\ q\equiv c_0 \, (\text{mod }c)}}\frac{1}{\phi(q)\phi(d)}\sum_{\substack{\psi \, (\text{mod }d)\\ \chi\, (\text{mod } q)\\\operatorname{cond}(\chi )\leq R }} \sum_{n\leq x}\Lambda(n)\chi(n\overline{a_1}a_2)\psi(n\overline{d_0})\\
& \quad\quad + \sum_{\substack{q\leq Q\\ (q,a_1a_2)=1\\ q\equiv c_0\, (\text{mod } c)}}\frac{1}{\phi(d)}\sum_{\psi \, (\text{mod } d)}\sum_{n\leq x}\Lambda(n)u_R(n\overline{a_1}a_2;q)\psi(n\overline{d_0}) :=\mathcal{S}_1+\mathcal{S}_2,
\end{align*}
say. We first show 
\begin{align}\label{s2bound}
\mathcal{S}_2\ll c x (\log x)^{O(1)}R^{-1/11}.
\end{align}
As explained in \cite[p.\ 719-720]{Dr}, after applying  Heath-Brown's  identity (most conveniently phrased in \cite[Lemma 5]{BFI}),  
it is  enough to prove 
\begin{equation}\label{HBbound}
\begin{split}
\frac{1}{\phi(d)}\sum_{\psi \, (\text{mod }d)}&\psi(\overline{d_0})\sum_{\substack{Q\leq q \leq 2Q\\ (q, a_1a_2)=1\\ q\equiv c_0 \, (\text{mod } c)}} \operatornamewithlimits{\sum\cdots \sum}_{\substack{(1-\Delta)M_i<m_{i}\leq  M_{i}\\ (1-\Delta)N_i<n_{i}\leq N_{i}\\1\leq i \leq j}}\Big(\prod_{i=1}^j\mu(m_{i})\psi(m_i)\Big)\\
& \times (\log n_1)\psi(n_1\cdots n_j) u_R(n_1m_1\cdots n_jm_j \overline{a_1}a_2;q) \ll cx(\log x)^{O(1)}R^{-1}
\end{split}
\end{equation}
where $j\in\{1, 2, 3, 4,5 \}$,   and $Q, M_i, N_i$ are real numbers such that $Q \leq x^{1/2 + \delta}$, $ \prod_{i}M_iN_i\ll x$,  $M_i\leq 2x^{1/5}$. We need the dissection in short intervals in order to approximate the hyperbolic size conditions and treat the boundary terms trivially. With $\Delta = R^{-1/11}$, we obtain \eqref{s2bound} from \eqref{HBbound} by an argument identical to 
 \cite[p.\ 720]{Dr}.  

The contribution from $\prod_i M_i N_i\leq x^{1-\delta}$ is trivially bounded by $O_\epsilon(x^{1-\delta+\epsilon})$. So we can assume that $\prod_{i}M_iN_i\geq x^{1-\delta}$. We rename $\prod_i M_i N_i=x$ for convenience and it is enough to prove \eqref{HBbound} for $Q\leq x^{1/2+4\delta},|a_2|\leq x^{1+4\delta}, |a_1|\leq x^{1+4\delta}$. Note that $M_i\leq x^{1/4-\delta}$ if $\delta$ is small enough. 
We write $M_i=x^{\mu_i}$ and $N_i=x^{\nu_i}$ and we have 
\begin{align}\label{partition}
 0\leq \mu_j< 1/4, \quad \mu_1+\cdots + \mu_j+ \nu_1+\cdots+ \nu_j=1.
\end{align} 
Let $\eta > 0$ be sufficiently small. 
If there exists an index $k$ such that $\mu_k$ or $\nu_k$ lies in the interval $[\eta, {1/4+\eta}]$, then we can apply Proposition \ref{TypeII} to see that 
\eqref{HBbound}  holds for $Q\leq x^{1+4\delta}, |a_2|\leq x^{4\delta}, |a_1|\leq x^{1+4\delta}$ as long as $\delta$ is small enough depending on $\eta$.  If none of the $\mu_i$ or $\nu_i$ lies in $[\eta, 1/4+\eta]$, then we combine all  $\mu_i$ and $\nu_i$ smaller than $\eta$ to  $\tau$ with $\tau< 10\eta$. This includes in particular $\mu_i$. The remaining $\nu_i$ must satisfy $\nu_i\geq 1/4+\eta$ and thus there are at most three of such terms. In summary, we see that \eqref{partition} can be partitioned as  
$\tau+ \nu_{i_1}+ \cdots +\nu_{i_s}=1,$ with $\tau<10\eta$, $\nu_{i_j}\geq 1/4+\eta$ and $s\leq 3$.  Since $(d, q)=1$, we see that \eqref{HBbound}  becomes 
\begin{align}
&\sum_{\substack{Q\leq q\leq 2Q\\ (q, a_1a_2)=1\\ q\equiv c_0 \, (\text{mod } c)}}\sum_{\substack{ln\asymp x\\ ln \equiv d_0\, (\text{mod }d)}}(\alpha_1*\cdots * \alpha_s)(n)\gamma(l)u_R(nl\overline{a_1}a_2;q)\label{d3case}\\
&=\sum_{\substack{Q\leq 2Q\\ (q, a_1a_2)=1\\ q\equiv c_0 \, (\text{mod } c)}}\Big(\hspace{-0.15cm}\sum_{\substack{ln \asymp x\\ ln \equiv d_0 \, (\text{mod } d)\\ ln \equiv a_1 \overline{a_2}\, (\text{mod } q)}}\hspace{-0.3cm}(\alpha_1* \cdots *\alpha_s) (n)\gamma(l) 
- \sum_{\substack{ln \asymp x}}\frac{(\alpha_1* \cdots *\alpha_s) (n)\gamma(l) }{\phi(q)\phi(d)}\sum_{\substack{\psi \, (\text{mod } d)\\ \chi \, (\text{mod } q)\\ \operatorname{cond}(\chi)\leq R}}\chi(nl\overline{a_1}a_2)\psi(nl \overline{d_0})\Big)\nonumber\\
&=\sum_{\substack{Q\leq q\leq 2Q\\ (q, a_1a_2)=1\\ q\equiv c_0 \, (\text{mod } c)}} \sum_{ln \asymp x}(\alpha_1* \cdots \alpha_s)\gamma(l)u_1(nl (\overline{a_1}a_2\bar{d}d+d_0\bar{q}q);qd)\nonumber\\&\quad\quad-\sum_{\substack{Q\leq q \leq 2Q\\ (q, a_1a_2)=1\\ q\equiv c_0 \, (\text{mod } c)}}\sum_{\substack{ln \asymp x}}\frac{(\alpha_1* \cdots *\alpha_s) (n)\gamma(l) }{\phi(q)\phi(d)}\sum_{\substack{\chi \,(\text{mod }q)\\ \operatorname{cond}(\chi)\leq R\\ \operatorname{cond}(\psi \chi)>1}}\chi(nl\overline{a_1}a_2)\psi(nl \overline{d_0})\nonumber,
\end{align}
where $\alpha_i$ is $\mathds{1}_{((1 - \Delta) N_i, N_i]}$ or $ \mathds{1}_{((1 - \Delta)N_i, N_i]}\log$, $s\leq 3$ and $\gamma(l)\ll \tau_{2j-1}(l)\log l$ with $l \leq x^{10\eta}$.  
For the first term we apply  \cite[Lemma 2]{BFI2} along with partial summation to remove the logarithm if necessary. This is a   deep input, since it uses Deligne's estimates for exponential sums over algebraic varieties over finite fields. In this way we obtain a power saving
$$\sum_{l \leq x^{10\eta}}\sum_{Q \leq q \leq 2Q}  \frac{x^{1-\delta_1} }{q} \ll x^{1-\delta_1 + 10\eta}$$
for some $\delta_1 > 0$ and $Q \leq x^{1/2 + \delta_1}$. 
The second term is 
$$
\ll \sum_{l\leq x^{10\eta}}\sum_{q\leq Q}\frac{|\gamma(l)|}{\phi(q)\phi(d)}\sum_{\psi \, (\text{mod } d)} \sum_{\substack{\chi \,(\text{mod }q)\\ \operatorname{cond}(\chi)\leq R\\ \operatorname{cond}(\psi \chi)>1}}\prod_{i=1}^s \left|\sum_{n_i\leq N_i}\chi(n_i)\psi(n_i)\alpha_i(n_i) \right |\ll x^{O(\eta)} Q^{1+\epsilon} (dR)^{3/2+\epsilon}.
$$
Combining  the two bounds, we have that 
\eqref{d3case} is bounded by $x^{1-\delta}$ for $Q\leq x^{1+4\delta}$ as long as $\eta$ sufficiently small and $\delta$ is small enough depending on $\eta$ and $\delta_1$.

Now we turn to $\mathcal S_1$. We extract the contribution where both $\chi$ and $\psi$ are trivial which gives the desired main term. For the remaining terms, we note that 
for any non-principal character $\chi_1$ (mod $q_1$) and any fixed positive constant $B$, we have uniformly for $q_1 \leq (\log x)^B$ that 
\begin{align*}
\sum_{n\leq x}\Lambda(n)\chi_1(n)\ll x e^{-c(B) \sqrt{\log x}},
\end{align*}
where $c(B)$ is some positive constant depending on $B$. 
Thus for $d\leq (\log x)^ C, R\leq (\log x)^B$, the remaining portion of $\mathcal{S}_1$ is 
\begin{align*}
\ll  \sum_{\substack{q\leq Q\\ (q, a_1a_2)=1\\ q\equiv c_0 \, (\text{mod } c)}}\frac{1}{\phi(qd)}\sum_{\psi \, (\text{mod } d)}\sum_{\substack{\chi \, (\text{mod } q)\\ \operatorname{cond}(\chi)\leq R\\\operatorname{cond}(\psi \chi)>1}}xe^{-c_2 \sqrt{\log x}}\ll  dRxe^{-c_3  \sqrt{\log x}} \ll xe^{-c_4  \sqrt{\log x}}
\end{align*}
where $c_2, c_3, c_4$ depend on $B, C$. Under GRH, we have for any non-principal character $\chi$ the bound
\begin{align*}
\sum_{n\leq x}\Lambda(n)\chi(n\overline{a_1}a_2)\ll x^{1/2+\epsilon}.
\end{align*}
Thus for $R\leq x^{\delta}$ the non-principal portion of $\mathcal{S}_1$ is 
$O(dRx^{1/2+\epsilon}). $ 
Combining with \eqref{s2bound}, we complete the proof.

\section{The Titchmarsh divisor problem and variations}\label{sec5}

We start with two elementary lemmas. The first is \cite[Lemma 13.1]{FGHM}.
\begin{lemma}\label{bfi-lemma} For $a\in \Bbb{Z} \setminus \{0\}$ and $x \geq 2$ we have
\begin{displaymath}
\begin{split}
&\sum_{\substack{m \leq x\\(m, a) = 1}} \frac{m}{\phi(m)} = x c_0(a) + O(\tau(a) \log x),\\
&\sum_{\substack{m \leq x\\(m, a) = 1}} \frac{1}{\phi(m)} =  c_0(a) ( \log x  + \gamma) + c_0'(a)+ O(\tau(a)x^{-1} \log x),
\end{split}
\end{displaymath}
where $c_s(a)$ is as in \eqref{c1c2}. 
\end{lemma}

\begin{lemma}\label{chi/phi}
	Let $\chi \pmod m$ be a non-principal Dirichlet character. Then for any non-zero integer $f$ and $x\geq 2$, we have
	\begin{align*}
	&\sum_{\substack{n\leq x\\ (n, f)=1}}\frac{\chi(n)}{\phi(n)}=L(1, \chi)c(\chi,f)+O\left( \frac{\tau(f)m^{1/2}\log m\log x}{x}\right),\\
	&\sum_{\substack{n\leq x\\ (n, f)=1}}\frac{\chi(n)n}{\phi(n)}\ll    \tau(f)m^{1/2}\log m\log x,
	\end{align*}
	where $c(\chi, f)$ is as in \eqref{c-chi}. 
\end{lemma}
 
\textbf{Proof.}  It is a simple consequence of the   P\'olya-Vinogradov inequality (and M\"obius inversion)  that
$$\sum_{\substack{n\leq x\\ (n, f)=1}}\frac{\chi(n)}{n}=\prod_{p\mid f}\Big(1- \frac{\chi(p)}{p}\Big)L(1, \chi)+O\left( \frac{\tau(f)m^{1/2}\log m}{x}\right).$$
Using the identity $n/\phi(n)=\sum_{d\mid n} \mu^2(d)/\phi(d)$ we conclude
	\begin{align*}
	\sum_{\substack{n\leq x\\ (n, f)=1}}{\frac{\chi(n)}{\phi(n)}}	&=\sum_{\substack{ d\leq x\\ (d, f)=1}}\frac{\mu^2(d)\chi(d)}{\phi(d)d}\sum_{\substack{n\leq x/d\\ (n, f)=1}}\frac{\chi(n)}{n}\\
	&=\sum_{\substack{ d\leq x\\ (d, f)=1}}\frac{\mu^2(d)\chi(d)}{\phi(d)d}\Big(\prod_{p\mid f}\Big(1-\frac{\chi(p)}{p}\Big)L(1, \chi)+O\Big(\frac{d}{x}\tau(f)m^{1/2}\log m \Big)\Big).
	\end{align*}
 We complete the sum over $x$ and obtain the first stated formula. The second one is similar (but easier).\\

We recall the Brun-Titchmarsh inequality \cite[Theorem 6.6]{IK}
\begin{align*}
	\sum_{\substack{x\leq n\leq x+y\\ n\equiv a\, (\text{mod } q)}}\Lambda(n)\ll \frac{y \log(x + y)}{\phi(q)\log (y/q)}.
	\end{align*}
	for $(a, q) = 1$, $q < y$. \\

We are now prepared for the \textbf{proof of Theorem \ref{thm-a}}. For notational simplicity we only display the case $\sigma=1$, the other case is essentially identical except for notational changes. 
As before let $R = (\log x)^{A}$ unconditionally or $R = x^{\delta}$ on GRH. 	The terms with $n < X (\log X)^{-2}R^{-1}$ can be estimated trivially, so that we can assume without loss of generality that $Y > X(\log X)^{-2}R^{-1}$. We choose a parameter $Q$ and write
\begin{displaymath}
\begin{split}
&\sum_{Y <n \leq X}  \Lambda(n) \tau(n+f)  =\sum_{Y <n \leq X}  \Lambda(n) \sum_{\substack{q \leq Q\\ q \mid n+f}} 1 + \sum_{Y <n \leq X}  \Lambda(n) \sum_{\substack{q < (n+f)/Q\\ q \mid n+f}} 1\\
& = \sum_{q \leq Q} \sum_{\substack{Y <n \leq X \\ n \equiv - f\, (\text{mod }q)}}  \Lambda(n)   +\sum_{q \leq (X+f)/Q}  \sum_{\substack{\max(Y, qQ - f) <n \leq X\\ n \equiv - f\, (\text{mod }q)}}  \Lambda(n)\\
&  = \sum_{q \leq Q} \sum_{\substack{Y <n \leq X \\ n \equiv - f\, (\text{mod }q)}}  \Lambda(n)   +\sum_{q \leq (Y+f)/Q}  \sum_{\substack{Y <n \leq X\\ n \equiv - f\, (\text{mod }q)}}  \Lambda(n) + 
\sum_{(Y+f)/Q < q \leq (X+f)/Q}  \sum_{\substack{ qQ - f <n \leq X\\ n \equiv - f\, (\text{mod }q)}}  \Lambda(n). 
 \end{split}
\end{displaymath}
If $Y + f < X^{1/2 + \delta/2}$ with $\delta$ as in Theorem \ref{PrimeAp}, we choose $Q = X^{1/2 + \delta/2}$. Otherwise we choose $Q = X^{1/2}$. 	The point of this maneuver is that either the condition $q \leq (Y+f)/Q$ is non-existent, or $(Y+f)/Q > X^{\delta/2}$, so that the error term in Lemma \ref{bfi-lemma} saves a power of $Q$. In order to apply Theorem \ref{PrimeAp} to the last term, we argue as in \cite[p.\ 720]{Dr} and split the $n$-sum in short subsums of the form $((1 - \Delta)N, N]$ with $\Delta = R^{-1/2}$. There are $O(\Delta^{-1} \log X)$ such intervals. The condition $qQ - f < n$ is only relevant in one such interval, whose contribution we estimate by the Brun-Titchmarsh inequality, getting a total error of $O(\Delta X (\log X)^2)$.  Otherwise we apply Theorem \ref{PrimeAp} with $c = d = a_2 = 1$ which gives a total error of $ \ll \Delta^{-1} X (\log X)^{1-A} \ll X (\log X)^{1 - A/2}$. It remains to assemble the main term, which by the prime number theorem up to an admissible error equals
\begin{displaymath}
\begin{split}
   \sum_{\substack{q \leq Q\\ (q, f) = 1}}  \frac{X-Y}{\phi(q)}   +\sum_{\substack{q \leq (Y+f)/Q\\ (q, f) = 1}} \frac{X-Y}{\phi(q)}  + 
\sum_{\substack{(Y+f)/Q < q \leq (X+f)/Q\\ (q, f) = 1}}  \frac{X - ( qQ - f )}{\phi(q)} .
 \end{split}
\end{displaymath}
The result follows now easily from Lemma \ref{bfi-lemma}, distinguishing the cases $(Y+f)/Q < 1$ and $(Y+f)/Q > 1$. \\


The \textbf{proof of Theorem \ref{hooley}} is similar, except that the derivation of the main term is slightly different. For $\sigma = 1$ we have 
\begin{displaymath}
\begin{split}
&\sum_{Y <n \leq X}  \Lambda(n) (\chi_1\ast \chi_2)(n+f) \\
& =\sum_{Y <n \leq X}  \Lambda(n) \sum_{\substack{q \leq  Q\\ q \mid n+f}} \chi_1(q)\chi_2\Big(\frac{n+f}{q}\Big) + \sum_{Y <n \leq X}  \Lambda(n) \sum_{\substack{q < (n+f)/Q\\ q \mid n+f}} \chi_2(q)\chi_1\Big(\frac{n+f}{q}\Big)\\
&  = \sum_{q \leq Q} \chi_1(q)\sum_{\substack{Y <n \leq X \\ n \equiv - f\, (\text{mod }q)}}  \Lambda(n) \chi_2\Big(\frac{n+f}{q}\Big)    +\sum_{q \leq (Y+f)/Q}   \chi_2  (q) \sum_{\substack{Y <n \leq X\\ n \equiv - f\, (\text{mod }q)}}  \Lambda(n) \chi_1\Big(\frac{n+f}{q}\Big) \\
&\quad\quad\quad + 
\sum_{(Y+f)/Q < q \leq (X+f)/Q}  \chi_2  (q)  \sum_{\substack{ qQ - f <n \leq X\\ n \equiv - f\, (\text{mod }q)}}  \Lambda(n) \chi_1\Big(\frac{n+f}{q}\Big).
 \end{split}
\end{displaymath}
We consider the first sum and write $q = q_1q_2$ with $(q_1, c_2) = 1$, $q_2 \mid c_2^{\infty}$ getting
\begin{displaymath}
\begin{split}
  &    \sum_{\substack{q_2 \mid c_2^{\infty}\\q_2 \leq Q}} \chi_1(q_2) \sum_{\substack{q_1 \leq Q/q_2\\ (q_1, c_2) = 1}} \chi_1\bar{\chi}_2(q_1) \sum_{b \, (\text{mod }c_2)} \chi_2(b)\sum_{\substack{Y <n \leq X \\ n \equiv - f\, (\text{mod }q_1q_2)\\ (n+f)/q_2 \equiv b\, (\text{mod }c_2)}}  \Lambda(n) \\
  & = \sum_{\substack{q_2 \mid c_2^{\infty}\\q_2 \leq Q}} \chi_1(q_2)\sum_{b \, (\text{mod }c_2)} \chi_2(b) \sum_{a \, (\text{mod }c_1c_2)} \chi_1\bar{\chi}_2(a)\sum_{\substack{q_1 \leq Q/q_2\\ (q_1, c_2) = 1\\ q_1 \equiv a \, (\text{mod }c_1c_2)}}  \sum_{\substack{Y <n \leq X \\ n \equiv - f\, (\text{mod }q_1)\\ n  \equiv bq_2 - f (\text{mod }q_2c_2)}}  \Lambda(n).  
   \end{split}
\end{displaymath}
Up to a negligible error, we can restrict to  
$(f, q_1q_2) = (bq_2- f, c_2)= 1$. Moreover, 
the contribution of $q_2 \geq K$ is trivially bounded by
$$\ll X (\log X) \sum_{\substack{q_2 \mid c_2^{\infty}\\ q_2 \geq K}} \frac{1}{q_2} \ll X (\log X) K^{-1/2} \tau(c_2),$$
so that we can restrict the size of $q_2$ to a large power of $\log X$ (or a small power of $X$ under GRH). As above we may assume that $Y$ is sufficiently large, we choose $Q \in \{X^{1/2}, X^{1/2 + \delta/2}\}$ as above, and we restrict the $n$-sum to small intervals $((1 - \Delta)N, N]$ to deal with the fact that the size condition for $n$ depends on $q$. Theorem \ref{PrimeAp} can now be applied. Using the prime number theorem, we obtain up to an admissible error term the  main term
\begin{displaymath}
\begin{split}
   & \sum_{\substack{q_2 \mid c_2^{\infty}\\q_2 \leq  K\\ (q_2, f) = 1}} \chi_1(q_2)\sum_{\substack{b \, (\text{mod }c_2)\\ (bq_2 - f, c_2) = 1}} \chi_2(b)  \sum_{\substack{q_1 \leq  Q/q_2\\ (q_1, f) = 1}}  \chi_1\bar{\chi}_2(q_1) \frac{X -Y}{\phi(q_1q_2c_2)} . 
      \end{split}
\end{displaymath}
Since $\chi_2$ is primitive, the $b$-sum is easy to evaluate (cf.\ \cite[(3.9)]{IK}):
\begin{align*}
\sum_{\substack{b\, (\text{mod } c_2)\\ (bq_2-f, c_2)=1}}\chi_2(b)=\sum_{d\mid c_2}\mu(d)\sum_{\substack{bq_2\equiv f\, (\text{mod } d)\\ b \, (\text{mod } c_2)}}\chi_2(b)=\left\{\begin{array}{ll}
\chi_2(f)\mu(c_2), & q_2=1,\\
0, &q_2>1,
\end{array}\right.
\end{align*}
so that our main term becomes
$$\frac{X-Y}{\phi(c_2)} \sum_{\substack{q_1 \leq  Q \\ (q_1, f) = 1}} \frac{ \chi_1\bar{\chi}_2(q_1)}{\phi(q_1)} ,$$
which we evaluate with Lemma \ref{chi/phi}. The other two terms as well as the case $\sigma = -1$ are treated analogously. This completes the proof of Theorem \ref{hooley}. \\

The \textbf{proof of Theorem \ref{thm-b}} follows from the special case $\sigma = -1$, $c_1 = 1$, $c_2 = 2$  by partial summation together with the well-known formula $L(1, \chi_{-4}) = \pi/4$. 

\section{Proof of Theorem \ref{lem5}}\label{sec6}

The first part of the argument is standard. Let us denote by $\mathcal{S}$ the left hand side of the claimed inequality.  Using \cite[(1.6)]{DI}, we have
\begin{displaymath}
\begin{split}
\mathcal{S} =\sum_{\substack{\gamma = s\sqrt{r} c  \\\text{for some } (c, r) = 1} }\frac{s\sqrt{r}}{\gamma} \sum_{m \asymp M} \alpha_m F\Big(m, \frac{\gamma}{s\sqrt{r}}\Big)S_{\infty, 1/s}(m, \pm n, \gamma)e\Big( \mp \frac{n\bar{s}}{r}\Big) .  
\end{split}
\end{displaymath}
We drop the exponential as $n, s, r$ are fixed. Write $P = CMZnrs$. We separate variables by Fourier inversion and write
$$F(m, c) = \int_{-\infty}^{\infty} G\Big(u, \frac{4\pi \sqrt{mn}}{s\sqrt{r} c}\Big) e(mu) du, \quad G(u, \xi) = \int_{-\infty}^{\infty} F\Big(t, \frac{4\pi \sqrt{tn}}{s\sqrt{r}\xi}\Big)e(-tu)dt.$$
The function $G$ is supported on $\xi \asymp X$ in the second variable, and by partial integration we see that
$$ \frac{d^{\nu}}{d\xi^{\nu}}G(u, \xi) \ll_{A, \nu}  M (X/Z)^{-\nu} (1 + |u|M/Z)^{-A}$$
for any $\nu, A\in \Bbb{N}_0$. We truncate the $u$-integral at $|u| \ll  P^{\varepsilon}Z/M$, pull the $u$-integration outside and estimate it trivially in the end. This will cost a factor $ZP^{\varepsilon}$, and we absorb the factor $e(mu)$ into the coefficient $\alpha_m$. From now on we replace $F(m, c)$ by $f( 4\pi \sqrt{mn}/\gamma)$ where $f$ is compactly supported on $[X, 2X]$ and satisfies $f^{(\nu)} \ll_{\nu} (X/Z)^{-\nu}$. 
Applying the Kuznetsov formula to the $\gamma$-sum along with standard bounds for the Bessel transforms (see e.g.\ \cite[Lemma 2.1]{P2}) we obtain a decomposition $\mathcal{S} = \mathcal{S}_{\text{Maa{\ss}}} + \mathcal{S}_{\text{hol}} +  \mathcal{S}_{\text{Eis}};$ here  
\begin{displaymath}
\begin{split}
\mathcal{S}_{\text{Maa{\ss}}} & \ll   s\sqrt{r}\sum_{m \asymp M} \alpha_m   \sum_{\substack{\phi \in \mathcal{B}(rs)\\ t_{\phi} \ll ZP^{\varepsilon} }} \check{f}^{\pm}(t_{\phi})  \frac{ |\rho_{\phi, \infty}(m) \rho_{\phi, 1/s}(\pm n)  | }{\cosh(\pi t_\phi)}
\end{split}
\end{displaymath}
where again the subscripts $\infty, 1/s$ denote the cusp and 
$$\check{f}^{\pm}(t) = \left\{\begin{array}{l} \displaystyle \frac{\pi}{\sinh(\pi t) } \int_0^{\infty} \frac{J_{2it}(x) - J_{-2it}(x)}{2i} f(x) \frac{dx}{x} \\ \displaystyle\frac{2}{\pi} \cosh(\pi t) \int_0^{\infty} K_{2it} f(x) \frac{dx}{x} \end{array}\right\} \ll_A P^{\varepsilon} (X/Z)^{-2|\Im t|}(1+ |t|/Z)^{-A}.$$
The terms $\mathcal{S}_{\text{hol}}$ and $\mathcal{S}_{\text{Eis}}$ are similar but easier, as the only difficulty is the treatment of the exceptional spectrum. 
By the Cauchy-Schwarz inequality we conclude
\begin{displaymath}
\begin{split}
\mathcal{S}_{\text{Maa{\ss}}} \ll P^{\varepsilon}  s\sqrt{r}&\Big( \sum_{\substack{\phi \in \mathcal{B}(rs)\\ t_{\phi} \ll ZP^{\varepsilon} }} \frac{1}{\cosh(\pi t_\phi)} \Big|\sum_{m \asymp M} \alpha_m\rho_{\phi, \infty}(m)\Big|^2\Big)^{1/2} \\
&\Big( \sum_{\substack{\phi \in \mathcal{B}(rs)\\ t_{\phi} \ll ZP^{\varepsilon} }}\frac{(X/Z)^{-4|\Im t_\phi|}}{\cosh(\pi t_\phi)}|\rho_{\phi, 1/s}(n)|^2 \Big)^{1/2}, 
\end{split}
\end{displaymath}
up to a negligible error from the truncation of $t_\phi$.  
Here we used the fact that regardless of whether $\phi$ is even or odd, we have
  $|\rho_{\phi, 1/s}(-n)|^2 = |\rho_{\phi, 1/s}(n)|^2$.  
We apply the large sieve \cite{DI} on the first factor getting 
$$\mathcal{S}_{\text{Maa{\ss}}}  \ll P^{\varepsilon} s\sqrt{r}\Big( \Big(Z^2 + \frac{M}{rs}\Big) \| \alpha \|^2\Big)^{1/2} \Big( \sum_{\substack{\phi \in \mathcal{B}(rs)\\ t_{\phi} \ll ZP^{\varepsilon} }}\frac{(X/Z)^{-4|\Im t_{\phi}|}}{\cosh(\pi t_{\phi})} |\rho_{\phi,1/s}(n)|^2 \Big)^{1/2}.$$
So far this is standard, but the interesting point is how to treat the second factor, for which Pitt \cite[Section 2]{P2} has a relatively elaborate argument. We proceed   differently and apply directly the Kuznetsov formula backwards (i.e.\ in the form of \cite[Theorem 16.3]{IK}). Here we can add the Eisenstein spectrum by positivity, and we choose a test function majorizing $(X/Z)^{-4|\Im t_{\phi}|}$ times the characteristic function on $t_{\phi} \ll ZP^{\varepsilon}$.  To this end we split the $t_{\phi}$-sum into several pieces. For the  treatment of the non-exceptional part we decompose the  $t_{\phi}$-sum into $O(P^{\varepsilon})$ regions $T-\Delta \leq| t_{\phi}| \leq T+\Delta$ with $\Delta = T^{1-\varepsilon}$ and use the standard test function
$$h(t) =  T^{-2}\Big(t^2 + \frac{1}{4}\Big)\Big[ \exp\Big( - \Big(\frac{t - T}{\Delta}\Big)^2 \Big) + \exp\Big( - \Big(\frac{t + T}{\Delta}\Big)^2 \Big)\Big], $$
see e.g.\ \cite[Section 3]{JM}. A detailed analysis  of the corresponding Bessel transform
\begin{equation}\label{kuz}
\frac{2i}{\pi} \int_{-\infty}^{\infty} J_{2it}(x) \frac{h(t) t}{\cosh(\pi t)} dt
\end{equation}
can be found in \cite[(3.12), (3.19)]{JM}. 
We only need to know that \eqref{kuz} is negligible for $x \ll T^{2} P^{-\varepsilon} $ and otherwise bounded by $T^{1+\varepsilon}$. 
In this way we obtain easily the bound 
$$\sum_{\substack{\phi \in \mathcal{B}(rs)\\ t_{\phi} \ll ZP^{\varepsilon}, t_{\phi} \in \Bbb{R} }}\frac{|\rho_{\phi,1/s}(n)|^2}{\cosh(\pi t_{\phi})}  \ll P^{\varepsilon}\Big(Z^2 + \max_{1 \leq T \leq ZP^{\varepsilon}} T\sum_{\gamma \ll nT^{-2}P^{\varepsilon}} \frac{|S_{1/s, 1/s}(n, n, \gamma)|}{\gamma}\Big). $$
Here $\gamma$ runs over the allowable moduli for Kloosterman sums to the pair of cusps $(1/s, 1/s)$ of level $rs$, and by   \cite[Lemma 2.5]{DI} these are the  multiples of $rs/(r, s) = rs$ since $(r, s) = 1$.

For the exceptional spectrum we choose a different test function, namely $$h(t) = ((X/Z)^{2it} + (X/Z)^{-2it})^2 (t^2 + (2A)^2)^{-2A}\prod_{k=0}^A \Big(t^2 + \frac{(2k+1)^2}{4}\Big)$$
for fixed, but large $A$ (which we think of as $1/\varepsilon$), and recall that $X/Z \ll 1$. 
This is a variation of  
 \cite[(16.57)]{IK}, and the product over $k$ is inserted to cancel the poles of $1/\cosh(\pi t)$ in a large horizontal strip. Note that this function is even, holomorphic in a wide horizontal strip,  majorizes $(X/Z)^{-4|\Im t_{\phi}|}$ on the exceptional spectrum and is non-negative on the regular spectrum. To bound \eqref{kuz} with this test function, we only need the crude, but uniform estimate \cite[8.411.8]{GR} 
 $$ J_{2it + c}(x) \ll_c e^{\pi |t|}  x^c , \quad t \in \Bbb{R},  x > 0, c \geq 0.$$
In this way we see (with $c=0$) that  the Kuznetsov transform \eqref{kuz} is uniformly bounded. On the other hand, for $x \ll P^{-\varepsilon} (X/Z)^2$ we can shift the $t$-contour down to $\Im t = -A$ without crossing poles to see that \eqref{kuz} is negligible for such $x$. Thus,  with this test function the Kuznetsov formula returns the bound
$$ \sum_{\substack{\phi \in \mathcal{B}(rs)\\ t_{\phi} \ll ZP^{\varepsilon}, t_{\phi} \in i\Bbb{R} }}\frac{(X/Z)^{-4|\Im t_{\phi}|}}{\cosh(\pi t_{\phi})} |\rho_{\phi,1/s}(n)|^2 \ll 1+ \sum_{\gamma \ll P^{\varepsilon}Z^2|n|/X^2} \frac{|S_{1/s, 1/s}(n, n, \gamma)|}{\gamma}$$
for the exceptional spectrum.  By Weil's bound \cite[Lemma 2.6]{DI} for these Kloosterman sums we finally get 
\begin{displaymath}
\begin{split}
\mathcal{S}_{\text{Maa{\ss}}}  &\ll  P^{\varepsilon} s\sqrt{r}\Big( \Big(Z^2 + \frac{M}{rs}\Big) \| \alpha \|^2\Big)^{1/2} \Big( Z^2 +  \frac{Z\sqrt{n(n, rs)}}{Xrs}\Big)^{1/2}. \\
 \end{split}
 \end{displaymath}
The estimation of $\mathcal{S}_{\text{hol}}$ and $\mathcal{S}_{\text{Eis}}$ is similar but much easier, because no exceptional eigenvalues need to be treated, and the theorem follows.

\section{Preparatory lemmas}\label{sec7}

\subsection{Exponential and character sums}
We start with a generalization of \cite[Theorem 3]{P1}.
\begin{lemma}\label{char} Let $\alpha, \beta, h, r, a\in \Bbb{Z}$, $c \in \Bbb{N}$. Then
$$\sum_{\substack{x, y, z\, (\text{{\rm mod }} c)\\ (xyz(z+h), c) = 1}} e\Big(\frac{a\bar{x} - r x - \bar{a} y + r y + \alpha x \bar{z} - \beta y\overline{(z+h)}}{c}\Big) \ll (\alpha - \beta, \alpha h, \beta h, c_1)^{1/2} c_2^2 c_1^{3/2 + \varepsilon} (a, c)$$
where $c = c_1c_2$ with $(c_1, c_2) = 1$, $c_1$ squarefree and $c_2$ power-full. 
\end{lemma}

\textbf{Proof.} The case $(a, c) = 1$ is \cite[Theorem 3]{P1}. By multiplicativity we may restrict our attention to prime power moduli $c = p^{\nu}$. If $\nu \geq 2$, we apply the same argument as in \cite[(1.1)]{P1}. For $\nu = 1$ and $p \mid a$, the character sum is trivially bounded by $O(p^2)$ unless $\alpha \equiv \beta \equiv a \equiv r \equiv 0$ (mod $p$), in which case we obviously cannot do better than a bound of size $O(p^3)$. \\

Next we state and prove an analogue of \cite[Proposition 1]{P1}:
\begin{lemma}\label{lem1}  Let $r, f, Z, X_1, X_2, X_3, G$ be as in  Theorem \ref{thm1}.  For $h, \alpha \in \Bbb{Z} \setminus \{0\}$, $c, \gamma \in \Bbb{N}$, $(c , \gamma) = 1$, $c\gamma  \gg \max(X_1, X_2, X_3)/Z$, we have
\begin{displaymath}
\begin{split}
&\sum_{k, l, m} G(k, l, m) e\Big( \frac{\alpha klm}{\gamma}\Big) S(rklm + f, h, c)\\
& \ll c^{\varepsilon}
 \Big(X_1X_2X_3Z\frac{(r, c)(c, f, h)^{1/2}}{c^{1/2}} \\
&\quad\quad  + (c, f)^{1/2} (Z^3\gamma^2c^{7/4} +  Z^{5/2}X_1^{1/2} (c\gamma)^{3/2} + Z^2(X_2X_3)^{1/2} \gamma c^{5/4} + ZX_1(X_2X_3)^{1/2} c^{1/2})c_{\square}^{1/4} \Big)\end{split}
\end{displaymath}
where $c_{\square}$ is the power-full part of $c$.
\end{lemma}

\textbf{Proof.} There are two differences to the statement of \cite[Proposition 1]{P1}: on the one hand we have an extra exponential $e(\alpha klm/\gamma)$, on the other hand we have a general number $f$ instead of $-1$. We can assume without loss of generality that $(\alpha, \gamma) = 1$, otherwise we cancel common factors. We follow the proof of \cite[Proposition 1]{P1}  and apply Poisson summation to
$$\Lambda^{\ast} = \sum_{\substack{k, l, m\\ (k, c) = 1}} G(k, l, m) e\Big( \frac{\alpha klm}{\gamma}\Big) S(dklm -a, b, c)$$
for non-zero $a, b$, $d\mid c$, $c \mid (c/d)^{\infty}$ in residue classes modulo $\gamma c$ getting
$$\frac{d}{c^2\gamma^2} \sum_{u, v, w} \widehat{G}\Big( \frac{u}{\gamma  c}, \frac{v}{\gamma c/d}, \frac{w}{\gamma c/d}\Big)   S_{\gamma}(u, v, w) \underset{x, y\, (\text{{\rm mod }} c)}{\left.\sum\right.^{\ast}} e\Big(\frac{a\bar{\gamma}^3\bar{x} - b\gamma^3x - u\bar{y}+dxyvw}{c}\Big) $$
where
$$S_{\gamma}(u, v, w) = \sum_{\substack{x, y \, (\text{mod } \gamma)\\ xy \equiv u \bar{c}^3 \bar{\alpha}\, (\text{mod } \gamma)}} e\Big( - \frac{vx + wy}{\gamma}\Big)\ll \gamma^{1+\varepsilon} .$$
Note that unlike in the proof of \cite[Proposition 1]{P1} we are not assuming $(ab, c) = 1$. The condition $(b, c) = 1$ is never used in \cite[Proposition  1]{P1}, while the lack of $(a, c) = 1$ is compensated by our more general result in Lemma \ref{char} that allows arbitrary $a$, not necessarily coprime with $c$. 
The bound for $\Lambda^{\ast}$ can now be achieved verbatim as on  \cite[pp.\ 396-398]{P1} keeping in mind that  $U, V, W$ are now longer by a factor $\gamma$ and that Lemma \ref{char} introduces an extra factor $(a, c)$ in the bound of $|\Lambda_3|^2$ on \cite[p.\ 398]{P1}.  The condition $c\gamma \gg (X_1, X_2, X_3)/Z$ is needed to ensure $U \gg 1$ which is used on \cite[p.\ 398]{P1}. This gives
\begin{displaymath}
\begin{split}
  \Lambda^{\ast}& \ll \Big(X_1X_2X_3Z\frac{d}{c^{1/2}} \\
  &+ (a, c)^{1/2} \Big(Z^3\frac{c^{7/4} \gamma^2}{d} + \frac{Z^{5/2}X_1^{1/2}(c \gamma)^{3/2} }{d} + Z^2(X_2X_3)^{1/2} c^{5/4} \gamma + ZX_1(X_2X_3c)^{1/2}\Big) c_{\square}^{1/4} \Big) c^{\varepsilon} .
\end{split}
\end{displaymath}
We now finish the proof as on \cite[p.\ 399]{P1}, keeping in mind that the Weil bound 
introduces an additional  factor 
$$(\gamma_1, f, h)^{1/2}(\gamma_3, f, h)^{1/2}   \leq (c, f, h)^{1/2}$$
which together with the extra factor $(a, c)^{1/2}$  in the bound of $\Lambda^{\ast}$ yields
$$(\gamma_1, f, h)^{1/2}(\gamma_3, f, h)^{1/2}  (\gamma_4, f)^{1/2} \leq (c, f)^{1/2}.$$
(In \cite{P1} we have $f=-1$, so these gcd's are trivial.) This completes the proof. \\

We recall the delta-symbol of Duke-Friedlander-Iwaniec, see e.g.\ \cite[Section 20.5]{IK}. 

\begin{lemma}\label{lem4} Let $T \geq 1$. There exists a function $\Delta_t(n)$ with the following properties: we have $\Delta_t(n) = 0$ unless $t \ll T^{1/2}$,  and for $ |n| \ll  T$ we have $\frac{d^j}{dn^{j}} \Delta_t(n) \ll_j T^{-1/2} |n|^{-j}$ for all $j\in \Bbb{N}_0$ and $$
\delta_{n=0} = \sum_{c \ll T^{1/2}} \frac{1}{c} \underset{a\, (\text{{\rm mod }} c)}{\left.\sum\right.^{\ast}}  e\Big( \frac{an}{c}\Big) \Delta_c(n).$$
\end{lemma}

\subsection{Automorphic forms}\label{sec72} As specified in the introduction, let $\phi$ be a (holomorphic or Maa{\ss}) cuspidal newform for $\Gamma_0(N)$. We denote by $\mu$ its archimedean conductor. 
We frequently use the Rankin-Selberg bound
$$\sum_{n \leq X} |\lambda(n)|^2 \ll (N\mu X)^{\varepsilon} X.$$
This allows us, for instance, to estimate
\begin{equation}\label{RS}
\sum_{n \leq X} |\lambda(n)|(n, c)^{1/2}  \ll (c N\mu X)^{\varepsilon} X
\end{equation}
for $c \in \Bbb{N}$.

Next we need the following Voronoi formula which is a slight re-statement of  \cite[Lemma 2.4]{BMN} in the case of trivial nebentypus, but arbitrary infinity type of $\phi$. For a Dirichlet character $\chi$ we generally write $c(\chi)$ for its conductor. 
\begin{lemma}\label{lem2} Let $F$ be a smooth function with compact support in $(0, \infty)$, $a \in \Bbb{Z}$, $q\in \Bbb{N}$ with $(a, q) = 1$. Write $N_1 = (N, q)$, $N_2 = N/N_1$, $q_2 = (N_2^{\infty}, q)$, $q_1 = q/q_2$, $$\frac{a}{q} = \frac{a_1}{q_1} + \frac{a_2}{q_2}.$$ Then
$$\sum_n \lambda(n) e\left(\frac{an}{q}\right) F(n) =  \sum_{\chi \, (\text{{\rm mod }} N_1)} \sum_{\substack{\rho \mid NN_1\\ (\rho, q_1) = 1}} c(a_2/q_2, \rho, \chi)\frac{1}{q_1}\sum_{\pm}  \sum_n  \lambda_{\chi}(n) e\Big(\mp \frac{\overline{a_1\rho} n}{q_1}\Big) \mathcal{F}^{\pm}\Big( \frac{n}{q_1^2 \rho}\Big)$$
where   $\lambda_{\chi}$ are the Hecke eigenvalues of the contragredient of the twist  $\phi \times \chi$, 
\begin{equation}\label{boundc}
c(a_2/q_2, \rho, \chi) \ll (qN)^{\varepsilon}  \frac{(Nc(\chi))^{1/2} }{\rho}
\end{equation}
and 
\begin{equation}\label{vor}
\mathcal{F}^{\pm}(y) = \int_0^{\infty} F(x) \mathcal{J }^{\pm}(4\pi \sqrt{xy} )dx
\end{equation}
with
$$\mathcal{J}^{+}(z) = \begin{cases} 2\pi i^k J_{k-1}(z), & \phi \text{ holomorphic},\\ \frac{\pi i}{\sinh(\pi t)} (J_{2it}(y) - J_{-2it}(y)), &\phi \text{ Maa{\ss}}, \end{cases}$$ 
$$\mathcal{J}^{-}(z) = \begin{cases} 0, & \phi \text{ holomorphic},\\ 4\cosh(\pi t) K_{2it}(y), &\phi \text{ Maa{\ss}}. \end{cases}$$ 
\end{lemma}

\emph{Remark:} If $q_2 = 1$, the formula can be greatly simplified: the $\chi$-sum can be dropped, the $\rho$-sum consists of a single term $\rho = N_2$, and the constant $c$ is $\pm N_2^{-1/2}$ with a sign depending only on $\phi$, see \cite[Appendix A]{KMV}. \\

\textbf{Proof.} The exact shape   of the formula and the bound \eqref{boundc} require some explanation. The sum over $\rho$ runs over $N_2q_2^2$, but since $q_2\mid N_1$, we relaxed it to $\rho \mid NN_1$, $(\rho, q_1) = 1$. It follows from \cite[Lemma 2.3]{BMN} in their notation that $C(f, a_2/q_2, m, \chi) \ll (qN)^{\varepsilon} c(\chi)^{-1/2}$, and the third last display in the proof of \cite[Lemma 2.4]{BMN} shows 
$$c(a_2/q_2, \rho, \chi)  \ll (q N)^{\varepsilon} \frac{\sqrt{N c(\chi)^2}}{c(\chi)^{1/2}\rho}. $$
The exact bound, however, is irrelevant for the purpose of this paper, as long as it is polynomial in $N$, we only include it for completeness. Finally we argue that $c(a_2/q_2, \rho, \chi)$ is independent of $a_1$. To this end, we first observe that $\chi$ runs only modulo $q_2$ (which was relaxed in \cite{BMN} to $N_1$). Next the proof of \cite[Lemma 2.4]{BMN} shows that in the case of  trivial central character, the dependence on $a_1$ is $\tilde{\chi}^2(a_1)$ where $\tilde{\chi}$ is the restriction of $\chi$ modulo a number divisible by $q_1$. Since $(q_1, q_2) = 1$, there is no dependence on $a_1$. \\

The relevant integral transforms are estimated in following lemma. 

\begin{lemma}\label{lem-new} Let $X > 0$, $Z \geq 1$ and let  $G$ be a smooth function with compact support in $[X, 2X]$ satisfying $G^{(\nu)} \ll_{\nu} (X/Z)^{-\nu}$ for all $\nu \in \Bbb{N}_0$. Then 
$$y^{\nu} \frac{d^{\nu}}{dy^{\nu}}\mathcal{G}^{\pm}(y) \ll_{A, \nu} X\Big(1 + \frac{1}{(Xy)^{2\Im \mu+\varepsilon }}\Big)\big( \mu (1 + \sqrt{Xy})\big)^{\nu}\Big(1 + \frac{yX}{Z^2} + \frac{yX}{\mu^2}\Big)^{-A}$$
for any $\nu, A \geq 0$. \end{lemma}

\textbf{Proof.} This follows by partial integration based on the formula \cite[8.472.3 with $m=1$]{GR} 
$$\int_0^{\infty} F(x) \mathcal{B}_{k}(\alpha \sqrt{x}) dx = \pm \frac{2}{\alpha} \int_0^{\infty}\Big(F'(x) \sqrt{x} - \frac{k}{2}\frac{F(x)}{\sqrt{x}}\Big)   \mathcal{B}_{k+1}(\alpha \sqrt{x}) dx$$
for $\alpha > 0$, $F\in C^{\infty}_0((0, \infty)$, $k\in \Bbb{C}$ and $\mathcal{B}_k \in \{J_k, K_k, Y_k\}$ and standard bounds for Bessel functions
$$J^{(\nu)}_{k-1}(x),  \quad \cosh(\pi t) K^{(\nu)}_{2it}(x), \quad \frac{J_{2it}^{(\nu)}(x)}{\cosh(\pi t)} \ll_{\nu} \mu^{\nu} (1 + 1/x)^{\nu + 2\Im \mu+\varepsilon}.$$
The derivatives are computed using \cite[8.486.11, 8.471.2]{GR} and then the bounds can be obtained from \cite[8.411.1,  8.432.5 with $z=1$, 8.411.13]{GR} (or alternatively from the Mellin integrals as in \cite[Proposition 9]{HM}).  \\

We obtain the following two corollaries for Hecke eigenvalues in residue classes. 
\begin{cor}\label{lem3} Let $X, Z\geq  1$, let $G$ be a smooth function supported on $[X, 2X]$ with $G^{(\nu)} \ll_{\nu} (X/Z)^{-\nu}$, $a, q \in \Bbb{N}$. Then
$$\sum_{n \equiv a\, (\text{{\rm mod }} q)} G(n) \lambda(n) \ll (Nq\mu Z)^{\varepsilon} N^{1/2 } (N, q)^{3/2}   ( Z^2 + \mu^2 ) q^{1/2 }.$$
\end{cor}

\textbf{Proof.} The sum in question, say $S$,  equals
$$S = \frac{1}{q} \sum_{rs = q} \underset{u\, (\text{mod } r)}{\left.\sum\right.^{\ast}}e\left(-\frac{ua}{r}\right) \sum_n e\left(\frac{un}{r}\right) G(n) \lambda(n).$$
We use the notation from Lemma \ref{lem2} and write
$$N_1 = (N, r), \quad N_2 = N/N_1, \quad r_2 = (N_2^{\infty}, r), \quad r_1 = r/r_2, \quad \frac{u}{r} = \frac{u_1}{r_1} + \frac{u_2}{r_2}.$$
By the Chinese Remainder Theorem we obtain the bound
$$S \ll \frac{(qN)^{\varepsilon}}{q} \sum_{rs = q} 
 \sum_{\chi \, (\text{{\rm mod }} N_1)} \sum_{\substack{\rho \mid NN_1\\ (\rho, r_1) = 1}} \frac{(Nc(\chi))^{1/2}}{r_1 \rho}\sum_{\pm}  \sum_n \Big| \lambda_{\chi}(n)  
 r_2S(\mp \bar{\rho} n, -a, r_1) \mathcal{G}^{\pm}\Big(\frac{n}{r_1^2 \rho}\Big)\Big|.
 $$
Using Lemma \ref{lem-new}, Weil's bound and \eqref{RS}, 
 we obtain 
 $$S \ll \frac{(qN\mu Z)^{\varepsilon}}{q} \sum_{rs = q} N_1 (NN_1)^{\varepsilon} (NN_1)^{1/2} r_2r_1^{3/2}(\mu^2 + Z^2)$$
 and the lemma follows. \\

\begin{cor}\label{lem3a} Let $r, f\in \Bbb{Z} \setminus \{0\}$, $X_1, Z \geq 1$.  Let $h$ be a smooth function with support in $[X_1, 2X_1]$ satisfying $h^{(\nu)} \ll_{\nu} (X_1/Z)^{-\nu}$ for all $\nu \in \Bbb{N}_0$. Write $|r|X_1 = X$.   Suppose that $h(n)$ vanishes unless  $rn + f$ is some dyadic interval $[\Xi, 2\Xi]$ with $X/Z \ll \Xi \ll X+|f|$. Then
$$ \sum_n h(n)  \lambda(rn + f)  \ll  (Nq\mu Z)^{\varepsilon} N^{1/2 } (N, r)^{3/2}    \Big( Z\Big(1 + \frac{|f|}{X}\Big) + \mu\Big)^2 |r|^{1/2 }.$$
\end{cor}

\textbf{Proof.} The sum in question equals
$$ \sum_{m \equiv f\, (\text{mod } |r|)} \lambda(m)h\Big( \frac{m-f}{r}\Big).$$
The function $\tilde{h} : m \mapsto h((m-f)/r)$ has support in $[\Xi, 2\Xi]$ and satisfies $\tilde{h}^{(\nu)} \ll (Z/X)^{\nu} = (Z\Xi/X)^{\nu} \Xi^{-\nu}$ and 
%
the statement follows   from Corollary \ref{lem3}.  
 

\section{Shifted convolution sums}\label{sec8}
In this section we prove Theorems \ref{thm1} and \ref{lem8} and some consequences. We continue to denote by $\phi$ a cusp form as at the beginning of Section \ref{sec72}. 

\subsection{Proof of Theorem \ref{thm1}}  We call $\mathscr{S}$ the sum on the left hand side of the claimed inequality.  
 Suppose without loss of generality $X_1 \geq X_2 \geq X_3$. Recall that $rklm + f$ is restricted to a dyadic  interval of size $\Xi$ where 
 \begin{equation}\label{xi}
   X/Z \ll \Xi \ll X + |f|.
   \end{equation} 
   On the one hand we can apply Corollary \ref{lem3a} getting
\begin{equation}\label{s1}
\begin{split}
\mathscr{S}& \ll (N\mu XZ |f|)^{\varepsilon} \sum_{l \asymp X_2} \sum_{m \asymp X_3}  N^{1/2+\varepsilon}(N, rlm)^{3/2}  \Big( Z\Big(1 + \frac{|f|}{X}\Big) + \mu\Big)^2 |rlm|^{1/2+\varepsilon} \\
&\ll  (N\mu XZ |f|)^{\varepsilon}  N (N, r)  \Big( Z\Big(1 + \frac{|f|}{X}\Big) + \mu\Big)^2 \frac{1}{|r|} \Big(\frac{X}{X_1}\Big)^{3/2} .\end{split}
\end{equation}
 On the other hand, we can apply Lemma \ref{lem4} with $T = \Xi$ getting
\begin{displaymath}
\begin{split}
\mathscr{S}
= \sum_{k, l, m}  & G(k, l, m) \sum_n\lambda(n) \sum_{c \ll \Xi^{1/2}} \frac{1}{c} \underset{a\, (\text{{\rm mod }} c)}{\left.\sum\right.^{\ast}}  e\Big( \frac{a(rklm+f-n)}{c}\Big)  H_{c, k, l, m}(n) 
\end{split}
\end{displaymath}
where
$$  H_{c, k, l, m}(n) = h(n)\Delta_c(rklm+f-n)$$
and $h$ is a smooth function   with support of size $\Xi$ and $h^{(\nu)} \ll_{\nu} \Xi^{-\nu}$. We apply Voronoi summation (Lemma \ref{lem2}) to the $n$-sum and use the notation from that lemma:
$$N_1 = (N, c), \quad N_2 = N/N_1, \quad c_2 = (N_2^{\infty}, c), \quad c_1 = c/c_2, \quad \frac{a}{c} = \frac{a_1}{c_1} + \frac{a_2}{c_2}.$$
In this way we obtain
\begin{displaymath}
\begin{split}
 & \sum_n\lambda(n)   \underset{a\, (\text{{\rm mod }} c)}{\left.\sum\right.^{\ast}}  e\Big( \frac{a(rklm+f-n)}{c}\Big)  H_{c, k, l, m}(n) \\
  &=    \underset{a_2\, (\text{{\rm mod }} c_2)}{\left.\sum\right.^{\ast}}    e\Big( \frac{a_2(rklm+f)}{c_2}\Big)  \sum_{\chi \, (\text{{\rm mod }} N_1)} \sum_{\substack{\rho \mid NN_1\\ (\rho, c_1) = 1}} c(-a_2/c_2, \rho, \chi)\\
  & \quad\quad\quad\quad\quad\quad\times \frac{1}{c_1}\sum_{\pm}  \sum_n  \lambda_{\chi}(n) S(rklm+f, \pm \bar{\rho} n, c_1)   \mathcal{H}_{c, k, l, m}^{\pm}\Big( \frac{n}{c_1^2 \rho}\Big),
\end{split}
\end{displaymath}
so that 
\begin{displaymath}
\begin{split}
 \mathscr{S} \ll& (\mu XN|f|)^{\varepsilon}\sum_{c \ll \Xi^{1/2}}  
  \sum_{\chi \, (\text{{\rm mod }} N_1)} \sum_{ \rho \mid NN_1 } \sum_n   \frac{|\lambda_{\chi}(n)| (N c(\chi))^{1/2}c_2}{cc_1\rho}\\
 &\max_{a_2\, (\text{mod }c_2)}\Big|\sum_{k, l, m} G(k, l, m) e\Big( \frac{a_2rklm}{c_2}\Big) S(rklm+f, \pm \bar{\rho} n, c_1)   \mathcal{H}_{c, k, l, m}^{\pm}\Big( \frac{n}{c_1^2 \rho}\Big)\Big|. 
  \end{split}
\end{displaymath}
We partition the sum into pieces where $c_1 \asymp C \ll \Xi^{1/2}$. 
By  Lemma \ref{lem-new} we see that
$$k^{\nu_1} l^{\nu_2} m^{\nu_3} \frac{\partial^{\nu_1}}{\partial k^{\nu_1}} \frac{\partial^{\nu_2}}{\partial l^{\nu_2}} \frac{\partial^{\nu_3}}{\partial m^{\nu_3}} \mathcal{H}_{c, k, l, m}^{\pm}(y) \ll_{\nu_1, \nu_2, \nu_3, A} \Xi^{1/2}\Big(1 + \frac{1}{(\Xi y)^{2\Im \mu +\varepsilon}}\Big) \Big(1  + \frac{y\Xi}{\mu^2}\Big)^{-A}$$
so that in particular we can truncate the $n$-sum at $n \ll X^{\varepsilon}C^2 \rho\mu^2/\Xi$ up to a negligible error. This implies in particular$$C \gg \frac{X^{-\varepsilon}\Xi^{1/2}}{\rho^{1/2}\mu} \gg \frac{X^{1/2 - \varepsilon}}{Z^{1/2} N\mu},$$
otherwise the $n$-sum is empty. Suppose that
\begin{equation}\label{suppose}
  X_1 \ll \frac{Z^{1/2}X^{1/2-\varepsilon}}{N\mu},
\end{equation}
so that  $C\gg X_1/Z$. Then we are   in a position to apply Lemma \ref{lem1} and after summing over $n$  (recall \eqref{RS}) we obtain
\begin{displaymath}
\begin{split}
\mathscr{S} &\ll  (\mu XN|f|)^{\varepsilon}\sum_{c \ll \Xi^{1/2}}      \sum_{ \rho \mid NN_1 }    \frac{ N^{1/2} N_1^{3/2}c_2}{cc_1\rho} c_1^2 \rho \mu^2  \Xi^{-1/2}\Big(X_1X_2X_3Z\frac{(r, c_1)}{c_1^{1/2}}\\
&\quad\quad  + (c_1, f)^{1/2} (Z^3c_2^2c_1^{7/4} +  Z^{5/2}X_1^{1/2} (c_2c_1)^{3/2} + Z^2(X_2X_3)^{1/2} c_2 c_1^{5/4} + ZX_1(X_2X_3)^{1/2} c_1^{1/2})(c_1)_{\square}^{1/4} \Big)\\
&\ll (\mu XN|f|)^{\varepsilon}\sum_{c \ll \Xi^{1/2}}         N^{1/2} N_1^{3/2}    \mu^2   \Xi^{-1/2}\Big(X_1X_2X_3Z\frac{(r, c_1)c_2^{1/2}}{c^{1/2} }\\
&\quad\quad  + (c_1, f)^{1/2}\Big(Z^3c_2^{1/4}c^{7/4} +  Z^{5/2}X_1^{1/2} c^{3/2} + Z^2(X_2X_3)^{1/2}  c^{5/4} + ZX_1(X_2X_3)^{1/2}  c^{1/2} \Big)(c_1)_{\square}^{1/4} \Big).
\end{split}
\end{displaymath}
Using Rankin's trick, it is not hard to see that for $\alpha > -1$, $0 \leq \beta \leq 3/4$ we have 
\begin{displaymath}
\begin{split}
\sum_{c \leq C}  c^{\alpha}c_{\square}^{1/4} (c, M^{\infty})^{\beta} &\leq \sum_{c }  c^{\alpha}c_{\square}^{1/4} (c, M^{\infty})^{\beta} \Big(\frac{C}{c}\Big)^{\alpha+1+\varepsilon}\\
& \ll  C^{\alpha+1+\varepsilon} \prod_{p\mid M} \sum_{k=0}^{\infty} \frac{1}{p^{k(1+\varepsilon - \beta) - \frac{1}{2} [\frac{k}{2}]}} \ll C^{\alpha+1} (CM)^{\varepsilon}.
\end{split}
\end{displaymath}
Using this for the last line of the previous bound for $\mathscr{S}$, we can execute the sum over $c$ (recall  $|r|X_1X_2X_3 = X$, $X_1 \geq X_2 \geq X_3$ and \eqref{xi}, so that $X_2X_3 \ll (X/|r|)^{2/3}$) getting
\begin{displaymath}
\begin{split}
\mathscr{S} & \ll  (\mu |f| NX)^{\varepsilon} ZN^2\mu^2 \Big( \frac{X/|r|}{\Xi^{1/4}} + Z^2 \Xi^{7/8} + Z^{3/2} X_1^{1/2} \Xi^{3/4} + Z\Big( \frac{X}{|r|}\Big)^{1/3} \Xi^{5/8}  +\frac{X_1^{1/2} X^{1/2}}{|r|^{1/2}} \Xi^{1/4}\Big)\\
& \ll (\mu | f| NX)^{\varepsilon}  ZN^2\mu^2 \Big( Z^2   (X + |f|)^{23/24} +  Z^{3/2}X_1^{1/2} (X + |f|)^{3/4}  
\Big).
\end{split}
\end{displaymath}
Combining this with \eqref{s1} (multiplying the fourth root of \eqref{s1} with the second term of the previous display raised to the power $3/4$), we obtain
\begin{displaymath}
\begin{split}
\mathscr{S} & \ll (\mu Z | f| NX)^{\varepsilon} (1 + |f|/X)^{5/48} Z^3N^2\mu^2     (X + |f|)^{23/24} .
 \end{split}
\end{displaymath}
This holds under the assumption \eqref{suppose}. If \eqref{suppose} does not hold, we apply \eqref{s1} alone getting
$$\mathscr{S} \ll N^{5/2} \mu^{3/2} \frac{(N, r)}{|r|} \Big(Z\Big( 1 + \frac{|f|}{X}\Big) + \mu\Big)^2 (X/Z)^{3/4}(\mu Z | f| NX)^{\varepsilon} .$$
Combining the two previous bounds completes the proof.

\subsection{Proof of Theorem \ref{lem8}}  We denote by $\mathcal{T}$ the left hand side of the claimed inequality. We write $m_1 = dm_3$, $m_2 = dm_4$ with $(m_3, m_4) = 1$. Note that $(D, m_3m_4) = 1$ since $m_1, m_2$ are squarefree. We write $N = N_1N_2$ where $N_2 = (N, (m_3m_4)^{\infty})$ and $N_1 = N/N_2$. 
Using again the squarefreeness of $m_1, m_2$, we have
\begin{displaymath}
\begin{split}
\mathcal{T} &= \sum_{\substack{k_1, k_2 \asymp \Xi\\ k_1 \equiv f \, (\text{mod }d)\\ m_4(k_1-f) = m_3(k_2 - f)}} \lambda(k_1) \lambda(k_2) g\Big( \frac{k_2-f}{m_4d}\Big)\\
&=\frac{1}{d} \sum_{\nu \, (\text{mod } d)}  \sum_{\substack{k_1, k_2 \asymp \Xi \\ m_4(k_1-f) = m_3(k_2 - f)}} e\Big(\frac{(k_1 - f)\nu}{d}\Big)\lambda(k_1) \lambda(k_2) g\Big( \frac{k_2-f}{m_4d}\Big) .
\end{split}
\end{displaymath}
We insert a redundant function $\tilde{g}(k_1)$ where $\tilde{g}$ has support on $\frac{1}{2}\Xi  \leq k_1 \leq 3\Xi$ and satisfies $\tilde{g}(k_1) = 1$ for $\Xi \leq k_1 \leq 2\Xi$ and $\tilde{g}^{(\nu)} \ll_{\nu} \Xi^{-\nu}$ for all $\nu \in \Bbb{N}_0$. 

We detect the condition $m_4(k_1-f) = m_3(k_2 - f)$ by Jutila's circle method, see e.g.\ \cite[Proposition 3.1]{Bl}.  Let $P = |f| X M Z N \mu$,  and let $Q = P^A$ for some (fixed but) large $A$ be a gigantic parameter, $\omega$ a smooth non-negative function with support in $[Q, 2Q]$ satisfying $\omega^{(\nu)} \ll Q^{-\nu}$, $\omega(q) = 1$ on $[4Q/3, 5Q/3]$ and let $\mathcal{Q}$ be the set of moduli $q \in [Q, 2Q]$ which are divisible by $ F:= [N_1,D^2]$ and coprime to $m_3m_4$. Note that this is possible, since $(DN_1,  m_3m_4) = 1$. We define 
$$\Lambda = \sum_{q\in \mathcal{Q}} \omega(q) \phi(q) = (NM)^{o(1)}Q^2/F$$
and $\delta = 1/Q$. Then $\mathcal{T}  = \mathcal{T}_1 + O(P^{O(1)}Q^{-1/2})$, cf.\ e.g. \cite[(4.3)]{Bl}, where
\begin{displaymath}
\begin{split}
 \mathcal{T}_1 =  \frac{1}{2\delta \Lambda} \sum_{q \in \mathscr{Q}}& \omega(q) \underset{a \, (\text{mod }q)}{\left.\sum\right.^{\ast}} \int_{-\delta}^{\delta} \frac{1}{d} \sum_{\nu \, (\text{mod } d)}  \sum_{k_1, k_2} e\Big(\frac{(k_1 - f)\nu}{d}\Big)\\
&e\Big((m_4(k_1-f) - m_3(k_2 - f))\Big(\frac{a}{q} + \eta\Big) \Big)\lambda(k_1) \lambda(k_2) \tilde{g}(k_1)g\Big( \frac{k_2-f}{m_4d}\Big)d\eta .
\end{split}
\end{displaymath}
We have $e( n\eta)= 1 + O(|n| \delta)  = 1 + O(|n|/Q)$ for $n \in \Bbb{R}$. Since $Q$ is very large, we can remove the $\eta$-integral entirely up to a negligible error and approximate $\mathcal{T}_1$  by  
\begin{displaymath}
\begin{split}
 \mathcal{T}_2  = \frac{1}{  \Lambda}& \sum_{q \in \mathscr{Q}} \omega(q) \underset{a \, (\text{mod }q)}{\left.\sum\right.^{\ast}}   \frac{1}{d} \sum_{\nu \, (\text{mod } d)}  e\Big(\frac{-f\nu}{d}+\frac{(m_3-m_4)fa}{q}\Big)\\
 &\sum_{k_1} e\Big( \frac{k_1\nu}{d} + \frac{m_4k_1a}{q} \Big) \lambda(k_1) \tilde{g}(k_1)   \sum_{k_2} e\Big( \frac{-m_3k_2a}{q}\Big) g\Big( \frac{k_2-f}{m_4d}\Big) . \end{split}
\end{displaymath}
We change variables $\nu \mapsto \nu m_4 a$ (note that $(m_4a, d) = 1$ and recall that $d \mid q$)  getting
\begin{displaymath}
\begin{split}
 \mathcal{T}_2  = \frac{1}{  \Lambda}& \sum_{q \in \mathscr{Q}} \omega(q) \underset{a \, (\text{mod }q)}{\left.\sum\right.^{\ast}}   \frac{1}{d} \sum_{\nu \, (\text{mod } d)}  e\Big(  \frac{(m_3-m_4(1 + \nu q/d))fa}{q}\Big)\\
 &\sum_{k_1} e\Big(  \frac{m_4k_1a (1 + \nu q/d)}{q} \Big) \lambda(k_1) \tilde{g}(k_1)   \sum_{k_2} e\Big( \frac{-m_3k_2a}{q}\Big) h_{f, m_4d}(k_2)
  \end{split}
\end{displaymath}
where $h_{f, m_4d}(k_2)
 = g((k_2-f)/(m_4d))$. Note that $h_{f, m_4 d}$ is supported on $[\Xi, 2\Xi]$ and satisfies
 \begin{equation}\label{h}
   h_{f, m_4 d}^{(\nu)} \ll_{\nu} \Big(\frac{Z}{MX}\Big)^{\nu} =\Big(\frac{Z\Xi}{MX}\Big)^{\nu} \Xi^{-\nu}
 \end{equation}
 for all $\nu \in \Bbb{N}_0$. 
Next we apply Voronoi summation (Lemma \ref{lem2}) in both the $k_1$ and the $k_2$ sum; note that $(1 + \nu q/d, q) = 1$ since $d^2 \mid D^2 \mid q$, and we recall that $N=N_1N_2$ with $N_1\mid q$, $(N_2, q) = 1$, so that the remark after Lemma \ref{lem2} applies. 
We have
\begin{displaymath}
\begin{split}
\sum_{k_1}(...) & = \frac{\kappa}{q\sqrt{N_2}} \sum_{\pm}\sum_{k_1} \lambda(k_1)e\Big( \mp  \frac{k_1\overline{N_2m_4 a} (1 - \nu q/d)}{q} \Big) \tilde{\mathcal{G}}^{\pm}\Big(\frac{k_1}{q^2N_2}\Big),\\
\sum_{k_2}(...) & = \frac{\kappa}{q\sqrt{N_2}} \sum_{\pm}\sum_{k_2} \lambda(k_2)e\Big( \pm  \frac{k_2\overline{N_2m_3 a }}{q} \Big)  \mathcal{H}_{f, m_4d}^{\pm}\Big(\frac{k_2}{q^2N_2}\Big),
\end{split}
\end{displaymath}
where $\tilde{\mathcal{G}}^{\pm}$ and $\mathcal{H}_{f, m_4d}$ are derived from $\tilde{g}$ and $h_{f, m_4d}$ by \eqref{vor} and $\kappa \in \{\pm 1\}$ depends only on $\phi$.  Here we used that $\overline{1 + \nu q/d} \equiv 1 - \nu q/d$ (mod $q$).  We now need to understand the exponential sum
$$  \frac{1}{d} \sum_{\nu \, (\text{mod } d)} \underset{a \, (\text{mod }q)}{\left.\sum\right.^{\ast}}   e\Big(  \frac{(m_3-m_4(1 + \nu q/d))f\overline{N_2}a}{q} - \sigma_1  \frac{k_1\overline{ m_4 a} (1 - \nu q/d)}{q}+\sigma_2\frac{k_2\overline{ m_3 a }}{q} \Big)  $$
for $\sigma_1, \sigma_2 \in \{\pm 1\}$. Verbatim the same argument as in \cite[Lemma 3.1]{P2} (which is the case $f\bar{N}_2 \equiv -1 \, (\text{mod } q)$) shows that this equals
$$S\big((m_3-m_4)f\overline{N_2}, \overline{m_3m_4}(\sigma_2 m_4k_2 - \sigma_1 m_3k_1), q\big) \delta_{\sigma_1m_3^2 k_1 \equiv \sigma_2 m_4^2 k_2 \, (\text{mod } D)}$$
for $D^2 \mid q$. We obtain
\begin{displaymath}
\begin{split}
 \mathcal{T}_2 = \frac{1}{ N_2 \Lambda}& \sum_{q \in \mathscr{Q}} \frac{\omega(q)}{q^2}  \sum_{\sigma_1, \sigma_2\in \{\pm 1\}} \sum_{\sigma_1m_3^2 k_1 \equiv \sigma_2 m_4^2 k_2 \, (\text{mod } D) } \lambda(k_1) \lambda(k_2)\\
 &S\big((m_3-m_4)f\overline{N_2}, \overline{m_3m_4}(\sigma_2 m_4k_2 - \sigma_1 m_3k_1), q\big)   \tilde{\mathcal{G}}^{\sigma_1}\Big(\frac{k_1}{q^2N_2}\Big) \mathcal{H}_{f, m_4d}^{\sigma_2}\Big(\frac{k_2}{q^2N_2}\Big). 
   \end{split}
\end{displaymath}
By Lemma \ref{lem-new} 
we may truncate the $k_1, k_2$ sum at (recall \eqref{h})
\begin{displaymath}
\begin{split}
& k_1 \leq \mathcal{K}_1 := P^{\varepsilon}\frac{Q^{2} N_2\mu^2}{\Xi}, \\
& k_2  \leq \mathcal{K}_2:=P^{\varepsilon}Q^{2} N_2\Big(\frac{\mu^2}{\Xi} + \frac{\Xi Z^2}{(MX)^2}\Big) \ll P^{\varepsilon}\frac{Q^2 N_2}{\Xi}\Big(\mu^2 + Z^2\Big(1 + \frac{|f|}{MX}\Big)^2\Big).
\end{split}
\end{displaymath}
By assumption we have $m_3 \not= m_4$. By standard bounds for the Ramanujan sum, the contribution of $\sigma_2 m_4k_2 = \sigma_1 m_3k_1$ is easily seen to be $O(P^{O(1)}Q^{-1})$ (cf.\ \cite[(4.3)]{Bl} or \cite[p.\ 754]{P2}) and hence negligible, so from now on we exclude these pairs $(k_1, k_2)$. Recalling the definition of $\mathscr{Q}$ and using a partition of unity for the $k_1, k_2$-sum, we are left with estimating 
\begin{displaymath}
\begin{split}
  \frac{1}{ N_2 \Lambda}&  \sum_{\sigma_1, \sigma_2\in \{\pm 1\}} \sum_{\rho\not= 0} \sum_{\substack{\sigma_1m_3^2 k_1 \equiv \sigma_2m_4^2 k_2 \, (\text{mod } D)\\  \sigma_2 m_4k_2 - \sigma_1 m_3k_1 = \rho }} \lambda(k_1) \lambda(k_2) \psi_1\Big(\frac{k_1}{K_1}\Big)\psi_2\Big(\frac{k_1}{K_2}\Big)\\
 &\sum_{\substack{(q, m_3m_4) = 1\\ F \mid q}} \frac{\omega(q)}{q^2} S\big((m_3-m_4)f , \overline{N_2m_3m_4}\rho, q\big) \tilde{\mathcal{G}}^{\sigma_1}\Big(\frac{k_1}{q^2N_2}\Big) \mathcal{H}_{f, m_4d}^{\sigma_2}\Big(\frac{k_2}{q^2N_2}\Big)   \end{split}
\end{displaymath}
with $K_1 \leq \mathcal{K}_1$, $K_2 \leq \mathcal{K}_2$ and smooth compactly supported functions $\psi_1, \psi_2$. 
Again by Lemma \ref{lem-new} we can separate variables by a Mellin inversion argument. 
Let 
\begin{displaymath}
\begin{split}
& \Psi_1(x) = \tilde{\mathcal{G}}^{\pm}(x) \psi_1\Big(\frac{x q^2N_2}{K_1}\Big) = \int_{(0)} \widehat{\Psi}_1(s)  x^{-s} \frac{ds}{2\pi i},\\
& \Psi_2(x) =  \mathcal{H}_{f, m_4d}^{\pm}(x) \psi_2\Big(\frac{xq^2N_2}{K_2}\Big) = \int_{(0)} \widehat{\Psi}_2(s)  x^{-s} \frac{ds}{2\pi i} ,\end{split}
\end{displaymath}
so that 
$$\psi_1\Big(\frac{k_1}{K_1}\Big)\psi_2\Big(\frac{k_1}{K_2}\Big)\tilde{\mathcal{G}}^{\sigma_1}\Big(\frac{k_1}{q^2N_2}\Big) \mathcal{H}_{f, m_4d}^{\sigma_2}\Big(\frac{k_2}{q^2N_2}\Big) =  \int_{(0)}\int_{(0)}\widehat{\Psi}_1(s_1) \widehat{\Psi}_2(s_2)  \frac{(q^2 N_2)^{s_1+s_2}} {k_1^{s_1} k_2^{s_2}} \frac{ds_1\, ds_2}{(2\pi i)^2}. $$
 Then Lemma \ref{lem-new} together with partial integration gives
\begin{displaymath}
\begin{split}
& \widehat{\Psi}_j(s) \ll_A  \Xi\Big(1 + \frac{1}{(\Xi K_j/q^2N_2)^{2\Im \mu +\varepsilon}}\Big)  \Big(1 + \frac{(1 + |s|)}{\mu(1 + \sqrt{\Xi K_j/q^2N_2})}\Big)^{-A}, \quad \Re s = 0
\end{split}
\end{displaymath}
for $j = 1, 2$ and $A\geq 0$ for the Mellin transforms. Estimating the Mellin integrals trivially, it will suffice to bound
\begin{displaymath}
\begin{split}
 \mathcal{T}_3 =& P^{\varepsilon} \Xi^2\prod_{j=1}^2\Big[ \mu\Big(\Big(\frac{\Xi K_j}{Q^2N_2}\Big)^{1/2} + \Big(\frac{\Xi K_j}{Q^2N_2}\Big)^{-2\Im \mu}\Big) \Big] \\
 & \sup_{|t_j| \ll W} \frac{1}{ N_2 \Lambda}  \Big| \sum_{\rho\not= 0} \alpha_{\rho}\sum_{ (c, m_3m_4) = 1 } \frac{\omega(Fc)}{(Fc)^2} S\big((m_3-m_4)f , \overline{N_2m_3m_4}\rho, Fc\big)  c^{2i(t_1+t_2) }\Big|\end{split}
\end{displaymath}
where 
$$W = P^{\varepsilon} \mu\Big(\mu+Z\Big(1 + \frac{|f|}{MX}\Big)\Big)$$
is an upper bound for $\mu(1 + \sqrt{\Xi K_j /q^2 N_2})$ and
$$\alpha_{\rho} = \sum_{\substack{\sigma_1m_3^2 k_1 \equiv \sigma_2m_4^2 k_2 \, (\text{mod } D)\\ K_1 \leq k_1 \leq 2K_1, K_2 \leq k_2 \leq 2 K_2\\ \sigma_2 m_4k_2 - \sigma_1 m_3k_1 = \rho }} \lambda(k_1) \lambda(k_2)  k_1^{-it_1} k_2^{-it_2}. $$
We finally apply a smooth partition of unity to the $\rho$-sum, so that $$|\rho| \asymp R\ll \frac{M(K_1+K_2)}{d}.$$  
We can quickly estimate the 2-norm of $\bm \alpha$:
\begin{displaymath}
\begin{split}
\| \bm \alpha \|^2 = \frac{1}{D} &\sum_{u_1, u_2 \, (\text{mod } D)}  \int_0^1 \sum_{k_1, k_3 \asymp K_1}  \sum_{k_2, k_4 \asymp K_2}  \lambda(k_1) \lambda(k_2) \lambda(k_3) \lambda(k_4)  \Big(\frac{k_3}{k_1}\Big)^{it_1} \Big(\frac{k_4}{k_2}\Big)^{it_2}\\
& e\Big( \frac{(\sigma_1m_3^2 k_1 -\sigma_2m_4^2 k_2)u_1}{D }+ \frac{(\sigma_1m_3^2 k_3 -\sigma_2m_4^2 k_4)u_2}{D}\Big)\\
&e\big((\sigma_2 m_4k_2 - \sigma_1 m_3k_1  - \sigma_2 m_4k_3 + \sigma_1 m_3k_4)\beta\big) d\beta.
\end{split} 
\end{displaymath}
By  a uniform Wilton bound \cite[Proposition 4 \& 5]{HM} and partial summation we have
$$\| \bm \alpha \|^2 \ll P^{\varepsilon}W^4 N^4 \mu^8  K_1K_2.$$
This is now in good shape for an application of Theorem \ref{lem5} with $$(r, s, n, m, M, C, Z) \leftarrow (N_2m_3m_4, F, (m_3- m_4)f, \rho, R, Q/F, W')$$
where $$W' = P^{\varepsilon}\Big(\mu + Z\Big( 1 + \frac{|f|}{MX}\Big)\Big)\Big(\mu + Z^{1/2}\Big( 1 + \frac{|f|}{MX}\Big)^{1/2}\Big)  \geq W $$ which ensures
\begin{displaymath}
\begin{split}
X &= \frac{(Mn)^{1/2}}{s \sqrt{r} C} = \frac{(R|(m_3-m_4)f|)^{1/2}}{(N_2m_3m_4)^{1/2} Q} \ll \frac{(|f|(K_1+K_2))^{1/2}}{N_2^{1/2} Q} \\
&\ll  \frac{|f|^{1/2}}{\Xi^{1/2}} \Big(\mu + Z\Big(1 + \frac{|f|}{MX}\Big)\Big) \ll \Big(1 + \frac{|f|}{MX}\Big)^{1/2} Z^{1/2} \Big(\mu + Z\Big( 1 + \frac{|f|}{MX}\Big)\Big) \ll W'
\end{split}
\end{displaymath}
 as required for the application of Theorem \ref{lem5}. 
It is now a matter of book-keeping. We obtain
\begin{displaymath}
\begin{split}
 \mathcal{T}_3 \ll & P^{\varepsilon} \Xi^2\mu^2\prod_{j=1}^2\Big(\Big(\frac{\Xi K_j}{Q^2N_2}\Big)^{1/2} + \Big(\frac{\Xi K_j}{Q^2N_2}\Big)^{-2\Im \mu}\Big) \frac{W^2N^2 \mu^4  (K_1K_2)^{1/2}}{ N_2 \Lambda FQ } \\
&W'\Big(W'(FR)^{1/2} +\frac{(W')^2 FN_2^{1/2}M}{d} + \frac{R^{1/4} ((m_3-m_4)f, N_2Fm_3m_4)^{1/4} (W'Q)^{1/2}}{N_2^{1/4} (M/d)^{1/2}} \\
&\quad\quad+\frac{(W')^{3/2} ((m_3-m_4)f, N_2Fm_3m_4)^{1/4} (QF)^{1/2}N_2^{1/4}(M/d)^{1/2}}{R^{1/4}}\Big).
 \end{split}
\end{displaymath}
Only the first and the third in the last parenthesis end up with a non-negative power of $Q$. These terms are clearly increasing in $K_1, K_2$ and $R$, so we obtain (up to a negligible error)
\begin{displaymath}
\begin{split}
 \mathcal{T}_3 \ll & P^{\varepsilon} \Xi^2  \mu^2W\frac{W^2N^2 \mu^4   W}{\Xi  Q } W' \Big(W' \frac{F^{1/2} M^{1/2}Q N_2^{1/2}W/\mu}{(\Xi d)^{1/2}}\\
& \quad\quad\quad +  \frac{M^{1/4}Q^{1/2} N_2^{1/4}(W/\mu)^{1/2}((m_3-m_4)f, N_2Fm_3m_4)^{1/4}(W' Q)^{1/2}}{(\Xi d)^{1/4}N_2^{1/4} (M/d)^{1/2}}\Big) \\
 \ll &  P^{\varepsilon} \Xi^{1/2} \mu^5 W^5(W')^2 N^2        \frac{F^{1/2} M^{1/2}  N_2^{1/2} }{  d^{1/2}}+  P^{\varepsilon} \Xi^{3/4}  \mu^{11/2}W^{9/2}(W')^{3/2}N^2     \frac{ ((m_3-m_4)f, N_2Fm_3m_4)^{1/4}  }{   (M/d)^{1/4}} \\
 \ll & P^{\varepsilon} \mu^{11/2} W^5(W')^2  N^{5/2} \Big( \frac{(MX + |f|)^{1/2}M^{1/2} D}{d^{1/2}} + \frac{(MX+|f|)^{3/4}D^{1/4} (f, Dm_3m_4)^{1/4} }{M^{1/4}}\Big) ,  \end{split}
\end{displaymath}
where we used
\begin{displaymath}
\begin{split}
((m_3-m_4)f, N_2Fm_3m_4) & \leq \frac{N}{d}((m_1-m_2)f, D^2m_3m_4)  = \frac{ND}{d}\Big( \frac{(m_1-m_2)}{(m_1 - m_2, d^{\infty})} f, Dm_3m_4\Big)\\
& = \frac{ND}{d} (f, Dm_3m_4).
\end{split}
\end{displaymath}
This completes the proof. 

\subsection{A bilinear estimate}
A Cauchy-Schwarz argument yields the following estimate for certain bilinear forms, an analogue of \cite[Theorem 1.2]{P2}.
\begin{prop}\label{lem9} Let $M, X, Z \geq 1$, $f\in \Bbb{Z} \setminus \{0\}$, $\sigma \in \{\pm 1\}$. Let $\alpha_m, \beta_n$ be bounded sequences supported on $[M, 2M]$, $[X, 2X]$ respectively. Suppose that $\alpha_m\beta_n$ vanishes unless $m$ is squarefree and $\sigma mn+f$ is in some dyadic interval $[\Xi, 2\Xi]$ with $MX/Z \ll \Xi \ll MX + |f|$. Then
\begin{displaymath}
\begin{split}
\mathscr{B} &:= \sum_{m} \sum_n \alpha_m \beta_n \lambda(\sigma mn + f) \\
&\ll (N\mu MX |f|)^{\varepsilon} \mu^{23/4}\Big(\mu + Z\Big(1 + \frac{|f|}{MX}\Big)\Big)^{4} N^{5/4}\\
&\quad\quad\quad\big( X(MX + |f|) + X    (MX + |f|)^{1/2}M^{5/2}  + X     (MX+|f|)^{3/4}     M^{7/4}\big)^{1/2}.
\end{split}
\end{displaymath}
\end{prop}

\textbf{Proof.}  By Cauchy-Schwarz we have
$$|\mathcal{B}|^2 \ll \|\bm \beta\|^2 \sum_n g(n) \Big|\sum_m \alpha_m \lambda(\sigma mn + f)\Big|^2 \ll X  \sum_{m_1,m_2} \alpha_{m_1} \overline{\alpha_{m_2}} \sum_n g(n) \lambda (\sigma m_1n + f)\lambda (\sigma m_2n + f)$$
for a suitable non-negative function $g$ with support in $[X, 2X]$ satisfying $g^{(\nu)} \ll_{\nu} (X/Z)^{-\nu}$ and such that $\sigma mn + f \in [\Xi, 2\Xi]$ whenever $\alpha_m g(n) \not= 0$. The diagonal term $m_1=m_2$ contributes at most
$$ \ll (N\mu MX |f|)^{\varepsilon} X(MX + |f|). $$
We define $d, D, m_3, m_4$ as in Theorem \ref{lem8} and estimate 
the remaining portion $m_1 \not = m_2$, say $\mathscr{D}$, using the bound of Theorem   \ref{lem8}. In this way we obtain
\begin{displaymath}
\begin{split}
\mathscr{D}\ll&\Psi X \sum_{d \ll M}  \sum_{\substack{d\mid D \mid d^{\infty}  \\ D \ll M}}  \sum_{\substack{m_3, m_4 \ll M/d\\ (m_3, m_4) = 1, m_3\not= m_4\\ (m_3 - m_4, d^{\infty}) = D/d}} \Big( \frac{(MX + |f|)^{1/2}M^{1/2} D}{d^{1/2}} + \frac{(MX+|f|)^{3/4}D^{1/4} (f, Dm_3m_4)^{1/4} }{M^{1/4}}\Big).
\end{split}
\end{displaymath}
where
$$\Psi = (|f|N\mu MX)^{\varepsilon}\mu^{23/2}  \Big (\mu + Z \Big(1 + \frac{|f|}{MX}\Big)\Big)^8 N^{5/2}. $$
The first term is at most
\begin{displaymath}
\begin{split}
& \ll \Psi X\sum_{d \ll M}  \sum_{\substack{d\mid D \mid d^{\infty} \\  D \ll M}}    \frac{(MX + |f|)^{1/2}M^{1/2} D}{d^{1/2}}  \frac{M^2}{dD}  \ll \Psi  X    (MX + |f|)^{1/2}M^{5/2} .
\end{split}
\end{displaymath}
For the second term we write $D' = D/d$ and observe
\begin{displaymath}
\begin{split}
\sum_{\substack{m_3, m_4 \ll M/d\\ m_3\not= m_4\\ m_3 \equiv m_4\,(\text{mod } D')}} (f, m_3m_4)^{1/4} \leq \sum_{g_1, g_2 \mid f} (g_1g_2)^{1/4} \sum_{\substack{m'_3  \ll M/dg_1, m'_4 \ll M/dg_2\\ g_1m'_3\not= g_2m'_4\\ g_1m'_3 \equiv g_2m'_4\,(\text{mod } D')}} 1.
\end{split}
\end{displaymath}
We can assume without loss of generality $D' \ll M/d$, otherwise the sum is empty.  We can also assume without loss of generality $g_2 \geq g_1$. We necessarily have $\frac{(g_1, D')}{(g_1, D', g_2)} \mid m_4'$ and the congruence determines  
$m_3'$ modulo $D'/(g_1, D')$. We can therefore bound the previous display by
$$\sum_{\substack{g_1, g_2 \mid f\\ g_2 \geq g_1}} (g_1g_2)^{1/4} \frac{M(g_1, g_2, D')}{dg_2(g_1, D')} \Big(\frac{M(g_1, D')}{dg_1 D'} + 1\Big) \ll \sum_{\substack{g_1, g_2 \mid f\\ g_2 \geq g_1}} (g_1g_2)^{1/4} \frac{M }{dg_2 } \frac{M }{d  D'}  \ll \frac{|f|^{\varepsilon}M^2}{dD} .$$
 Thus the second term is bounded by
 \begin{displaymath}
\begin{split}
  \ll& \Psi X \sum_{d \ll M}  \sum_{\substack{d\mid D \mid d^{\infty}  \\ D \ll M}}    \frac{(MX+|f|)^{3/4}D^{1/4} (f, D)^{1/4} }{M^{1/4} }\frac{M^2}{dD} \ll  \Psi X     (MX+|f|)^{3/4}     M^{7/4}.
\end{split}
\end{displaymath}
Combining the previous estimates completes the proof.\\

\emph{Remark:} For $|f| \ll MX$ this is of the same quality as \cite[Theorem 1.2]{P2}. However, we only claim it for squarefree $m$, and we do not know  how to generalize this to non-squarefree $m$, as the argument in \cite[Section 6]{P2} uses Deligne's bound in a serious way. Luckily our application only requires squarefree $m$. Noting that the trivial bound is $\mathscr{B}  \ll (N\mu MX |f|)^{\varepsilon} (MX+|f|)$, we can drop the middle term,  and the bound of the lemma can be slightly weakened and simplified as
\begin{displaymath}
\begin{split}
\mathscr{B}  \ll (N\mu MX |f|)^{\varepsilon} \mu^{10}Z^{4} N^2\big(1 + \frac{|f|}{MX}\Big)^{5}\Big(M^{1/2} X  +   M^{5/4} X^{7/8} \big). 
\end{split}
\end{displaymath}
Up to unimportant factors of $\mu, N, Z, 1 + |f|/MX$, this is slightly stronger than \cite[(1.11)]{P2}.\\

\begin{cor}\label{cor10}  Let $M_1, M_2, X, Z \geq 1$, $f\in \Bbb{Z} \setminus \{0\}$, $\sigma \in \{\pm 1\}$. Let $\alpha_m, \beta_n$ be bounded sequences and suppose that $\alpha_m$ is supported only on squarefree integers. Let $g$  be a function with support in $[X, 2X]$ and $g^{(\nu)} \ll_{\nu}  (X/Z)^{-\nu}$ for all $\nu \in \Bbb{N}_0$. Suppose that $g(x)$ vanishes unless $\sigma x+f$ is in some dyadic interval $[\Xi, 2\Xi]$ with $X/Z \ll\Xi \ll X + |f|$.  
Then
\begin{displaymath}
\begin{split}
&\sum_{M_1 \leq m \leq M_2} \sum_{n} \alpha_m \beta_n g(mn) \lambda(\sigma mn+f) \ll (N\mu X |f|)^{\varepsilon}\mu^{10}Z^{5} N^2\Big(1 + \frac{|f|}{X}\Big)^{5}\Big(M_1^{-1/2} X  +   M_2^{3/8} X^{7/8} \Big).
\end{split}
\end{displaymath}
\end{cor}

This is the analogue of \cite[Corollary 1.3]{P2} and follows from Proposition \ref{lem9} simply by separating variables via 
$$g(nm) = \int_{(0)} \widehat{g}(s) (nm)^{-s} \frac{ds}{2\pi i}$$
and estimating the integral over $s$ trivially (which costs a factor $Z$). 

\section{The endgame}\label{sec9}
  
Verbatim the same as in \cite[Section 9]{P2}, based on Vaughan's identities and Theorem \ref{thm1}, Corollary \ref{lem3a} and Corollary \ref{cor10} (which is only applied for squarefree $m$) we obtain

  \begin{prop}\label{main2} Let $X, Z > 1$, $f\in \Bbb{Z} \setminus \{0\}$, $\sigma \in \{\pm 1\}$, $g$ a smooth function with support in $[X, 2X]$ such that $g^{(\nu)} \ll (X/Z)^{-\nu}$. Suppose that $g(x)$ vanishes unless $\sigma x+f$ is in some dyadic interval $[\Xi, 2\Xi]$ with $X/Z \ll\Xi \ll X + |f|$. Then there exist $B, \delta > 0$ such that 
$$\sum_{n} g(n) \Lambda(n) \lambda(\sigma n + f) \ll (\mu N Z (1+ |f|/X))^{B}X^{1 - \delta}.$$
\end{prop}

In order to obtain Theorem \ref{thm-c}, we use \eqref{RS} in the form
$$\sum_{N_1 < n \leq N_2} |\lambda(n)| \ll N_2^{1/2+\varepsilon} (N_2 - N_1)^{1/2}$$
We estimate the portion $n \leq X^{1-\delta_0}$ trivially, and for the rest we apply a smooth dyadic partition of unity where each constituent has support on some interval $[\Xi, 2\Xi]$ and is 1 on $[\Xi (1+1/Z), 2(1-1/Z)\Xi]$ for some very small power $Z$ of $X$.

\end{document}